\documentclass[11pt]{article}

\title{\vskip-1.0em\sc Dual convolution for the affine group of the real line}

\author{\sc Y. Choi, M. Ghandehari}
\date{3rd March 2021}

\usepackage[a4paper,left=27mm,right=27mm,top=30mm,bottom=30mm]{geometry}

\usepackage{sectsty}
\allsectionsfont{\sffamily}

\usepackage[colorlinks=true,citecolor=red,linkcolor=blue]{hyperref}

\input{fourier-affine_mac}

\begin{document}
\maketitle

\begin{abstract}
The Fourier algebra of the affine group of the real line has a natural identification, as a Banach space, with the space of trace-class operators on $L^2({\mathbb R}^\times, dt/ |t|)$.
In this paper we study the ``dual convolution product'' of trace-class operators that corresponds to pointwise product in the Fourier algebra. Answering a question raised in work of Eymard and Terp, we provide an intrinsic description of this operation which does not rely on the identification with the Fourier algebra, and obtain a similar result for the connected component of this affine group.
In both cases we construct explicit derivations on the corresponding Banach algebras, verifying the derivation identity directly without requiring the inverse Fourier transform. We also initiate the study of the analogous Banach algebra structure for trace-class operators on $L^p({\mathbb R}^\times, dt/ |t|)$ for $p\in (1,2)\cup(2,\infty)$.


%43A40 (1973-now) Character groups and dual objects 
%46J (1973-now) Commutative Banach algebras and commutative topological algebras
%46J99 (1973-now) None of the above, but in this section
%47B90 (2020-now) Operator theory and harmonic analysis 
\end{abstract}

\begin{section}{Introduction}

\begin{subsection}{Background and motivation}
Given a bounded and SOT-continuous representation $\pi$ of a topological group $G$ on a Banach space~$E$, one may associate to each $\xi\in E$ and $\phi\in E^*$ the {\itshape coefficient function} $\phi( \pi(\underline{\quad})\xi) \in C_b(G)$. The vector space generated by all coefficient functions of $\pi$ admits a natural norm, stronger than the uniform norm of $C_b(G)$, and its completion in this norm is called the {\itshape coefficient space of~$\pi$}. 

If $G$ is locally compact, we denote by $\FA(G)$ the coefficient space of the left regular representation $\lambda:G\to \cU(L^2(G))$. Eymard \cite{eymard64} showed that $\FA(G)$ is actually a Banach algebra with respect to pointwise product, now called the {\itshape Fourier algebra} of~$G$.
When $G$ is abelian, the Fourier transform gives an isometric isomorphism between $\FA(G)$ and the convolution algebra $L^1(\widehat{G})$.
Even when $G$ is non-abelian, a well-established theme in abstract harmonic analysis has been to view $\FA(G)$ as some kind of convolution algebra on a ``quantum group'' that is dual to~$G$. However, in most cases this ``dual convolution'' is only defined in a formal or abstract sense.

This article studies a particular case where this notion of dual convolution can be made precise and described explicitly.
Consider the group of affine transformations of $\Real$, given the natural topology, which we denote by $\affR$. This group has an unusual property that never occurs for non-trivial compact or abelian groups: writing $\sH=L^2(\Rx, dt/|t|)$, there is an \emph{irreducible} unitary representation $\pi:\affR\to\cU(\sH)$ such that $\FA(\affR)$ coincides with the coefficient space of~$\pi$, which we denote by $\FA_\pi(\affR)$.
Associated to $\pi$ is a surjective norm-decreasing map $\Psi:\sH\ptp\sHbar\to\FA_\pi(\affR)$, which is isometric since $\pi$ is irreducible.

Since $\FA_\pi(\affR)=\FA(\affR)$, and since $\FA(\affR)$ is a Banach algebra with respect to pointwise product, we can use the surjective isometry $\Psi:\sH\ptp\sHbar\to \FA_\pi(\affR)$ to equip $\cS_1(\sH)=\sH\ptp\sHbar$ with a commutative Banach algebra structure.
In \cite[Probl\`eme 2.7]{eymard-terp}, after making this observation, Eymard and Terp pose the following challenge:
\begin{quote}
``\textit{Interpr\'eter cette multiplication en terme des op\'erateurs!}''
\end{quote}
The present paper answers their challenge by providing an \emph{explicit formula} for the new multiplication on $\cS_1(\sH)$
 --- this is what we refer to as ``dual convolution'' for $\affR$.
To our knowledge, such a formula has not been recorded before in the literature.

Having established this explicit formula, the rest of our article investigates some applications and variations, described in more detail in Section~\ref{ss:outline}. These applications and variations are intended to demonstrate that the resulting Banach algebra $\cA$ can be studied directly, without any prior knowledge of the isomorphism $\Psi:\cA\to\FA(\affR)$, and to argue that $\cA$ is an object of intrinsic interest.
A loose but instructive parallel is with certain naturally occuring Banach function algebras, such as ${\rm AC}([0,1])$, that can be modelled as $L^1$-convolution algebras of certain semigroups.

Informally: by introducing dual convolution on $\cS_1(\sH)$, we are swapping an object where the algebra structure is easy to describe but the norm is complicated, for one where norm estimates are straightforward but the algebra structure is more complicated.
This offers an alternative point of view on $\FA(\affR)$, which could shed new light on its known properties as a Banach algebra. Moreover, analogous constructions for higher-dimensional semidirect product groups may yield new results for their Fourier algebras.
\end{subsection}

\begin{subsection}{Outline of our paper}\label{ss:outline}

Section~\ref{s:prelim} sets up the basic notation and definitions that will be used throughout the paper. We give an explicit definition/description of the group $\affR$ and the key representation $\pi:\affR\to\cU(\sH)$, and collect some known facts from the literature for ease of reference.

In Section~\ref{s:fusion-DC} we give an explicit formula for dual convolution as a bilinear map $\boxtimes: \cS_1(\sH)\times\cS_1(\sH)\to\cS_1(\sH)$. The formula is motivated by showing how one expresses the product of two coefficient functions of $\pi$ as a continuous average of other coefficient functions (a so-called ``fusion formula''). We show by explicit calculations, without invoking the representation~$\pi$, that $\boxtimes$ is commutative and associative. We also show that if trace-class operators on $\sH$ are given as integral kernel functions, then $\boxtimes$ can be described on that level also.

Writing $\cA$ for the Banach algebra $(\cS_1(\sH),\boxtimes)$: in Section~\ref{s:derivation} we construct a derivation $D:\cA\to\cA^*$ which has interesting operator-theoretic properties as a linear map between Banach spaces (it is cyclic, weakly compact, and ``co-completely bounded'' in the terminology of \cite{choi_Z2-pre}).
Usually, in constructing derivations on function algebras, it is easy to see that the derivation identity holds on a dense subalgebra, but hard to show that one has a well-defined and bounded map on the whole algebra. By working with dual convolution on $\cA$, the situation is reversed: it is easy to check that $D$ is a bounded linear map with the extra properties mentioned above, and the hard part is to verify the derivation identity.

The group $\affR$ is not connected, but has an index $2$ subgroup isomorphic to the semidirect product $\affRpos$, which is a fundamental example of a non-unimodular connected Lie group. (The notation will be explained in Section~\ref{s:prelim}.)
Since $\FA(\affRpos)$ cannot be identified with the coefficient space of a single irreducible representation, a direct description of dual convolution for $\affRpos$ is less straightforward.
In Section~\ref{s:subalgebra} we identify an explicit subalgebra of $\cA$ that corresponds to $\FA(\affRpos)$, and hence obtain an analogue of dual convolution for~$\affR$.
 We then show how the construction in Section~\ref{s:derivation} yields a derivation on the Fourier algebra of $\FA(\affRpos)$, which offers a new perspective on some resuts in~\cite{Choi-Ghandehari14}.

In Section~\ref{s:Lp-version} we consider $\FA^p_\pi$, the coefficient space of the $L^p$-analogue of $\pi$, from the viewpoint of dual convolution. We sketch how our explicit formula for $\boxtimes$ may be extended from $\cS_1(\sH)$ to $\cS_1(L^p(\Rx))$ for $1<p<\infty$, making $\cS_1(L^p(\Rx))= L^p(\Rx)\ptp L^q(\Rx)$ into a commutative Banach algebra~$\cA_p$.
We then show that $\FA^p_\pi$ is a Banach algebra in its natural norm and is isomorphic to $\cA_p$ (Theorem~\ref{t:Psi_p is injective}). Perhaps surprisingly, for $p\neq 2$ there is a crucial difference from the $p=2$ case: $\FA_\pi^p$ is not the same as the $L^p$-version of the Fourier algebra (Theorem~\ref{t:not A_p}), and it appears to be a new Banach function algebra about which we know little at this stage.

Finally, in Section~\ref{s:conclusion}, we make some remarks about possible directions for future work, and pose some explicit questions about the algebra~$\FA_\pi^p$.
In the appendix we show how the tensor product of two induced representations may be expressed as a direct integral of a family of induced representations, and use it to give an alternative proof of the fusion formula for coefficient functions of~$\pi$.
\end{subsection}
\end{section}

\begin{section}{Preliminaries}\label{s:prelim}

\begin{subsection}{Notation and some general background}
\label{ss:general}
If $\cH_1$ and $\cH_2$ are Hilbert spaces then $\cH_1\tp^2 \cH_2$ denotes their Hilbert-space tensor product.

Given a complex vector space $V$, the conjugate vector space $\overline{V}$ is defined to have the same underlying additive group as $V$, equipped with the new $\Cplx$-action $c\odot \xi = \overline{c}\xi$.
Note that if $\cH$ is a Hilbert space then the function $\cH\times\overline{\cH}\to\Cplx$ defined by $(\xi,\eta)\mapsto \ip{\xi}{\eta}$ is \emph{bilinear} rather than sesquilinear.

The symbol $\ptp$ denotes the projective tensor product of Banach spaces.
If $\cH$ is a Hilbert space then there is a standard identification of $\cH\ptp\overline{\cH}$ with the space $\cS_1(\cH)$ of trace-class operators on~$\cH$, defined by viewing the elementary tensor $\xi\tp\eta$ as the rank-one operator $\alpha \mapsto \ip{\alpha}{\eta}\xi$; this correspondence is an isometric, $\Cplx$-linear isomorphism of Banach spaces.

Coefficient functions associated to continuous bounded group representations were already defined in the introduction, but we did not give a precise definition of the corresponding coefficient spaces. Most of this article concerns unitary representations on Hilbert spaces, so we review some standard material here in order to fix our notation.

If $\sigma:G \to \cU(\cH)$ is a continuous unitary representation and $\xi,\eta\in \cH$, we denote the associated {\itshape coefficient function}
 $x\mapsto \ip{\sigma(x)\xi}{\eta}$ by  $\xi*_\sigma\eta \in C_b(G)$.
There is a contractive, \emph{linear} map $\Psi_\sigma: \cH\ptp\overline{\cH} \to C_b(G)$ defined by $\Psi_\sigma(\xi\tp\eta) = \xi*_\sigma\eta$.
We denote the range of $\Psi_\sigma$ by $\FA_\sigma(G)$, or simply $\FA_\sigma$ if the group $G$ is clear from context; this is the {\itshape coefficient space of~$\sigma$}, and we equip it with the quotient norm pushed forward from $\cH\ptp\overline{\cH} /\ker(\Psi_\sigma)$.
%% It can be shown that with this norm, $\Psi_\sigma)^*$ maps $\FA_\sigma(G)^*$ isometrically onto the von Neumann algebra generated by $\pi(G)$ inside $\Bdd(\cH)$.

Two special cases should be singled out:
\begin{enumerate}
\item
If $\lambda$ denotes the left regular representation $G\to \cU(L^2(G))$, then $\FA_\lambda(G)$ coincides with the {\itshape Fourier algebra} of $G$, and is usually denoted by $\FA(G)$. (This is not Eymard's original definition of $\FA(G)$ but the equivalence is proved in \cite[Ch.~3]{eymard64}; see also \cite[Prop.~2.3.3]{KL_AGBGbook}.)
With our definition, the fact that $\FA(G)$ is closed under pointwise product follows from Fell's absorption principle.
\item
If $\sigma:G \to \cU(\cH)$ is \emph{irreducible}, then $\Psi_\sigma:\cH\ptp\overline{\cH}\to \FA_\sigma(G)$ is injective, hence is an isometric isomorphism of Banach spaces.
This result is due to Arsac; the proof combines a duality argument (see e.g.~\cite[Lemma 2.8.2]{KL_AGBGbook}) with Schur's lemma for irreducible \emph{unitary} representations.
\end{enumerate}
Moreover, if $\sigma'$ is a direct sum of countably many copies of $\sigma$, then $\FA_{\sigma'}(G)=\FA_\sigma(G)$. (See e.g.~\cite[Prop.~2.8.8]{KL_AGBGbook}.)

\begin{rem}
The space $\FA_\sigma(G)$ was originally introduced by Arsac but defined in a different way, as the closed linear span of $\{\xi*_\sigma\eta \colon \xi,\eta\in \cH\}$ inside the {\itshape Fourier--Stieltjes algebra} ${\rm B}(G)$. We will not discuss ${\rm B}(G)$ in this paper; the equivalence of this original definition with our one can be found in e.g.~\cite[Theorem 2.8.4]{KL_AGBGbook}.  
\end{rem}

\end{subsection}

\begin{subsection}{The affine group of $\Real$}
\label{ss:affine-group}
$\Rx$ denotes the multiplicative group of $\Real$, equipped with the subspace topology; it has a Haar measure $dt/|t|$ where $dt$ denotes usual Lebesgue measure on~$\Real$.
We write $\Rxpos$ for the subgroup of $\Rx$ consisting of strictly positive real numbers; the notation is consistent with using $G_e$ to denote the connected component of a locally compact group~$G$.

When dealing with $L^p$-spaces on $\Rx$, we will usually omit mention of the Haar measure and merely write $L^p(\Rx)$; this should not be confused with $L^p(\Real)$ which always means the $L^p$-space for the Lebesgue measure on~$\Real$.

We define $\affR$ to be the set $\{(b,a)\colon b\in \Real, a\in\Rx\}$ equipped with the product topology of $\Real\times\Rx$ and the following multiplication:
\begin{equation}
(b,a)\cdot (b',a') = (ab'+b,aa')
\end{equation}
With this choice, the map $(b,a) \mapsto \twomat{a}{b}{0}{1}$ is a homomorphism $\affR \to \GL_2(\Real)$.
Inversion in $\affR$ is given by
\begin{equation}
(b,a)^{-1} = ( -b/a, 1/a)
\end{equation}
Note that $\Real$ embeds as a normal closed subgroup of $\affR$ via $b\mapsto (b,1)$, while $\Rx$ embeds as a closed subgroup via $a\mapsto (0,a)$.

In harmonic analysis it is more common to work with the subgroup $\{(b,a)\colon b\in\Real, a\in\Rxpos\}$. This is a connected Lie group, often referred to in the literature as ``the real $ax+b$ group''; we shall return to it in Section~\ref{s:subalgebra}.
\end{subsection}

\begin{subsection}{The key representation and its coefficient space}

As in the introduction, we let $\sH$ denote $L^2(\Rx)$. There is a continuous unitary representation $\Pi: \affR \to \cU(\sH)$, defined by
\begin{equation}\label{eq:induced}
\Pi(b,a)\xi(t) \defeq e^{2\pi i bt^{-1}}\xi(a^{-1}t) \qquad(\text{$b\in\Real$, $a\in\Rx$; $\xi\in\sH$, $t\in \Rx$}).
\end{equation}
This is a special case of a more general construction: if we consider the character $\chi_1$ on $\Real$ given by
$\chi_1(t) = \exp(2\pi it)$, the previous formula may be written as
\begin{equation}
\Pi = \Ind_{\Real}^{\affR} \chi_{1}\;,
\end{equation}
where we use the explicit realization of an induced representation for a semidirect product group, as described in ``Realization III'' of \cite[section 2.4]{KTbook} (see Appendix \ref{app:induced} for details). 
Mackey theory tells us that $\Pi$ is irreducible, and is the only infinite dimensional irreducible representation of $\affR$.

In this article we work not with $\Pi$ but with a unitarily equivalent form (which matches the representation defined in \cite[Equation~(1.3)]{eymard-terp}).
For a $\Cplx$-valued function on a group $G$, define $\check{f}:G\to\Cplx$ by $\check{f}(x)=f(x^{-1})$.
Since Haar measure on $\Rx$ is invariant under the change of variables $t\leftrightarrow t^{-1}$,
the map $\xi\mapsto\check{\xi}$ defines an isometric involution $J:\sH\to\sH$.
 We now define $\pi =  J\Pi(\cdot)J :\affR\to\cU(\sH)$. Explicitly, given $\xi\in\sH$ and $b\in\Real$, $a\in\Rx$, we have
\begin{equation}
\pi(b,a)\xi(t) \defeq e^{2\pi i bt}\xi(ta) \quad\qquad(t\in \Rx).
\end{equation}

We claimed in the introduction that $\FA_\pi(\affR)=\FA(\affR)$.
This can be seen as follows.
%%%%%%%%%%%%%%%%%%%%%%%%%%%%%%%
%c.f.~the remarks at the end of \cite[Ch.~1, \S1]{eymard-terp}.
%
%By continuity of induction (see e.g.~\cite[Proposition 5.42]{KTbook}), the singleton $\{\pi\}$ in $\widehat{G}$ is dense in the Fell topology.
%%%%%%%%%%%%%%%%%%%%%%%%%%%%%%
%
The left regular representation $\lambda$ of $\affR$ can be obtained by inducing the left regular representation of ${\mathbb R}$, which we denote by $\lambda_{\Real}$.
Note that $\lambda_{\Real}$ is unitarily equivalent to a direct integral (over $\Rx$) of all nontrivial characters of $\Real$. Moreover, each such character is induced to a representation of $\affR$ equivalent to $\pi$.
Since induction and direct integration commute, it follows that
$\lambda$
is equivalent to $\pi\tp I_{\cH}$ for some separable Hilbert space~$\cH$.
% can be written as a direct sum of countably many copies of $\pi$.
 Hence $\pi$ is weakly equivalent with $\lambda$,
and $\FA_\pi(\affR)=\FA_\lm(\affR)=\FA(\affR)$ by the results mentioned in Section~\ref{ss:general}.

\begin{rem}
The equality $\FA_\pi(\affR)=\FA(\affR)$ implies that $\FA_\pi(\affR)$ is closed under pointwise product. In Section~\ref{s:fusion-DC} we will give an alternative proof of this fact, using dual convolution on $\sH\ptp\sHbar$. In Section~\ref{s:Lp-version} we will see that this alternative proof carries over to the $L^p$-analogue of $\FA_\pi(\affR)$, but that this space is \emph{not} equal to the $L^p$-analogue of $\FA(\affR)$.
\end{rem}

We shall write $\Psi$ rather than $\Psi_\pi$ for the canonical quotient map $\sH\ptp\sHbar\to\FA_\pi(\affR)$, $\xi\tp\eta\mapsto \xi*_\pi\eta$.
Since $\pi$ is irreducible, $\Psi$ is injective by the remarks in Section~\ref{ss:general}, although we shall not use this fact when defining dual convolution in Section~\ref{s:fusion-DC}.

\begin{rem}
\begin{romnum}
\item
In \cite{eymard-terp}, the map $\Psi$ is denoted by $\overline{\cF}$ and called ``{\itshape la co-transformation de Fourier}'' for the group $\affR$. Note that because $\affR$ is non-unimodular, composing $\Psi$ with the operator-valued Fourier transform $\cF:f \mapsto \pi(f) = \int_{\affR} f(x) \pi(x) \,dx$ does not yield the map $f\mapsto \check{f}$, and so (as observed in \cite{eymard-terp}) $\Psi$ should not be called an ``inverse Fourier transform''. However, the philosophy of Fourier inversion guides much of what we do in this article.

\item
For most of our article, the fact that $\FA_\pi(\affR)=\FA(\affR)$ does not play a big role in our calculations, since we are not relying on the modified Plancherel formula for this group.
The exceptions are in Section~\ref{s:subalgebra}, where we use general facts about Fourier algebras of open subgroups,
and in the proof of Theorem~\ref{t:not A_p}, where we use results of Herz on Fig\`a-Talamanca--Herz algebras.

\end{romnum}
\end{rem}
\end{subsection}

\begin{subsection}{Bochner integrals and related measure theory}
\label{ss:bochner background}
Our explicit formula for dual convolution is expressed as a Bochner integral, which requires attention to questions of strong measurability (also referred to in the literature as {\itshape Bochner measurability}).
A very thorough treatment of strong measurability and the Bochner integral can be found in \cite[Section~1.2.b]{HVVW_vol1}.

It is usually impractical to verify directly that a given Banach-space valued function is strongly measurable.
For functions with values in an $L^p$-space an alternative approach is provided by the following result: given two sigma-finite measure spaces $(\Om_1,\mu_1)$ and $(\Om_2,\mu_2)$, and $1\leq p < \infty$, there is a natural embedding
\[
L^p(\Om_1,\mu_1)\tp L^p(\Om_2,\mu_2) \longrightarrow L^p(\Om_1,\mu_1 ; L^p(\Om_2,\mu_2))
\]
where $f\tp g$ is sent to the function $\omega_1\mapsto f(\omega_1)g$.
This embedding extends to an isometric isomorphism of Banach spaces $L^p(\Om_1\times\Om_2,\mu_1\times\mu_2) \cong L^p(\Om_1, \mu_1; L^p(\Om_2,\mu_2))$ (see e.g.\ \cite[Prop.\ 1.2.24]{HVVW_vol1} for the proof of a more general statement).
In particular, elements of $L^p(\Om_1\times \Om_2,\mu_1\times\mu_2)$ define  strongly $\mu_1$-measurable functions $\Om_1\to L^p(\Om_2,\mu_2)$.

\end{subsection}
\end{section}

%%%%%%%%%%%%%%%%%%%%%%%%%%%%%%%%%%%%%%%%%%%%%%%

\begin{section}{Fusion and dual convolution}\label{s:fusion-DC}
\paragraph{Notation.}
For $r\in\Rx$, let $\lm(r):\sH\to \sH$ denote the usual ``left translation'' by~$r$ (multiplicative in this context), i.e. $\lm(r)\xi(t) = \xi(r^{-1}t)$.
Similarly $\rho(r):\sH\to\sH$ denotes ``right translation'' by $r^{-1}$, i.e. $\rho(r)\xi(t)=\xi(tr)$.

We use both $\lm$ and $\rho$, even though $\Rx$ is abelian, because we have in mind possible extensions of the following calculations to semidirect products of the form $\Real^n\rtimes D$ where $D\subset\GL_n(\Real)$ need not be abelian.

\begin{subsection}{An explicit formula for fusion of coefficients}
To avoid any doubt we shall pay close attention to issues of convergence and integrability.

Let $\xi_1,\eta_1,\xi_2,\eta_2\in\sH$. For each $(b,a)\in\affR$,
\[
\ip{\pi(b,a)\xi_1}{\eta_1}
\ip{\pi(b,a)\xi_2}{\eta_2}
=
\int_{\Real} e^{2\pi i bt} \xi_1(ta)\overline{\eta_1(t)} \frac{dt}{|t|}
\int_{\Real} e^{2\pi i bs} \xi_2(sa)\overline{\eta_2(s)}\,\frac{ds}{|s|}
\]
(where as usual we treat a measurable function defined on $\Rx$ as a measurable function defined on $\Real$, by prescribing some arbitrary value at $0$).

Let $d(t,s)$ denote the Haar measure on $\Real^2$.
Observe that the function 
\[
(t,s)\mapsto \frac{e^{2\pi i bt}}{\abs{t}\abs{s}} \xi_1(ta)\overline{\eta_1(t)} 
e^{2\pi i bs} \xi_2(sa)\overline{\eta_2(s)}
\]
is integrable on $\Real^2$, since by Tonelli's theorem for $\Real^2$ followed by Cauchy--Schwarz for $\sH$,
\[
\begin{aligned}
\int_{\Real^2}  \abs{ \xi_1(ta)} \abs{ \eta_1(t)} 
\abs{ \xi_2(sa) } \abs{ \eta_2(s)}  \, \frac{d(t,s)}{\abs{t}\abs{s}}  
 & =
 \int_{\Real}  \abs{ \xi_1(ta)} \abs{ \eta_1(t)} \, \frac{dt}{|t|} 
\;\int_{\Real} \abs{ \xi_2(sa) } \abs{ \eta_2(s)} \,   \frac{ds}{|s|} \\
& \leq
\norm{\rho(a)\xi_1}_{\sH}\norm{\eta_1}_{\sH}
\norm{\rho(a)\xi_2}_{\sH}\norm{\eta_2}_{\sH}
<\infty.
\end{aligned}
\]
Therefore, the following changes of variable and order of integration are valid:
\begin{equation}\label{eq:getting started}
\begin{aligned}
& %\phantom{=}
\ip{\pi(b,a)\xi_1}{\eta_1}
\ip{\pi(b,a)\xi_2}{\eta_2}
\\
\expo{Fubini}\quad & =
\int_{\Real^2} e^{2\pi i bt} \xi_1(ta)\overline{\eta_1(t)}
\, e^{2\pi i bs} \xi_2(sa)\overline{\eta_2(s)}\,\frac{d(t,s)}{\abs{t}\abs{s}}
\\
 \expo{$t\mapsto t-s$}\quad &= 
\int_{\Real^2} e^{2\pi i b(t-s)} \xi_1((t-s)a)\overline{\eta_1(t-s)} \,
 e^{2\pi i bs} \xi_2(sa)\overline{\eta_2(s)}\,\frac{d(t,s)}{\abs{t-s}\abs{s}}
 \\
 \expo{$s\mapsto tu$} \quad & = 
\int_{\Real^2} e^{2\pi i bt}
 \xi_1((1-u)ta) \xi_2(uta) \, \overline{\eta_1((1-u)t)\eta_2(ut)}\,\frac{dt\,du}{\abs{t}\abs{1-u} \abs{u}}
\\
 \expo{Fubini} \quad & = 
\int_{-\infty}^\infty \left( \int_{\Rx} e^{2\pi i bt}
 \xi_1((1-u)ta) \xi_2(uta) \, \overline{\eta_1((1-u)t)\eta_2(ut)}\,\frac{dt}{|t|}\right) \frac{du}{\abs{1-u} \abs{u}}
\\
\end{aligned}
\end{equation}

One can now show that for fixed $u\in \Real\setminus \{0,1\}$, the inner integral in the last line of Equation~\eqref{eq:getting started} can be written as $\ip{\pi(b,a)\alpha_u}{\beta_u}$ for suitable $\alpha_u,\beta_u\in \sH$, and that
\[
\int_{\Real\setminus\{0,1\}} \norm{\alpha_u}_{\sH} \norm{\beta_u}_{\sH} \,\frac{du}{\abs{1-u}\abs{u}}
\leq
\norm{\xi_1}_{\sH}
\norm{\xi_2}_{\sH}
\norm{\eta_1}_{\sH}
\norm{\eta_2}_{\sH}
\]
so that $(\xi_1*_\pi\eta_1)\cdot(\xi_2*_\pi \eta_2)$ is a weighted average of explicit coefficient functions $\alpha_u*_\pi \beta_u$ as $u$ varies; this is what we mean by a ``fusion formula'' for coefficient functions.

For technical reasons, we first make a further change of variables $u\mapsto 1- (1+h)^{-1}$.
Then $|1-u|^{-1}|u|^{-1}du = |h|^{-1}dh$, and so the last line of Equation~\eqref{eq:getting started} is equal to
\begin{equation}\label{eq:fusion step1}
%\begin{aligned}
%\ip{\pi(b,a)\xi_1}{\eta_1}
%\ip{\pi(b,a)\xi_2}{\eta_2} \\
% = 
%
\int_{-\infty}^\infty \left( \int_{\Rx} e^{2\pi i bt}
\xi_1\bigl(\frac{ta}{1+h}\bigr) \xi_2\bigl(\frac{ta}{1+h^{-1}}\bigr) \,
 \overline{\eta_1\bigl(\frac{t}{1+h}\bigr) \eta_2\bigl(\frac{t}{1+h^{-1}}\bigr)}
\,\frac{dt}{|t|}\right) \,\frac{dh}{|h|}
%
% \end{aligned}
\end{equation}

\begin{lem}\label{l:fund unitary}
Given $X\in \sH\tp^2\sH = L^2(\Rx\times\Rx)$, define $V(X): \Rx\times\Rx \to\Cplx$ by
\[
V(X)(h,t) = X\left( \frac{t}{1+h} , \frac{t}{1+h^{-1}}\right).
\]
Then $V(X)\in L^2(\Rx\times\Rx)$, and $V:L^2(\Rx\times\Rx) \to L^2(\Rx\times\Rx)$ is an isometry.
\end{lem}

\begin{proof}
Clearly $V(X)$ is measurable. Then
\[ \begin{aligned}
& \phantom{\quad} \int_{\Rx\times\Rx}
\abs{V(X)(h,t)}^2 \, \frac{d(h,t)}{\abs{h}\abs{t}} 
\\
& = \int_{\Rx\times\Rx} 
\left\vert X\left( \frac{t}{1+h} , \frac{t}{1+h^{-1}} \right) \right\vert^2 \, \frac{d(t,h)}{\abs{t}\abs{h}}  
\\
& = \int_{\Rx} \left( \int_{\Rx}
\abs{ X( t, th ) }^2  \, \frac{dt}{|t|} \right)\, \frac{dh}{|h|}
  & \expo{Tonelli, then $t\mapsto t(1+h)$ }
  \\
& = \int_{\Rx} \left( \int_{\Rx}
\abs{ X( t, h ) }^2  \, \frac{dh}{|h|} \right)\, \frac{dt}{|t|}
  & \expo{Tonelli, then $h\mapsto t^{-1}h$ }
  \\
& = \int_{\Rx\times\Rx}
\abs{ X( t, h ) }^2  \, \frac{d(h,t)}{\abs{h}\abs{t}} \;.
  &\expo{Tonelli}
\end{aligned}
\]
Thus $V(X)\in L^2(\Rx\times\Rx)$ and $V$ is an isometry, as required.
\end{proof}

Note that if $X\in C_c(\Rx\times\Rx)$ then so is $V(X)$. However, if $X\in C_c(\Rx)\tp C_c(\Rx)$, we see no reason to expect that $V(X)\in C_c(\Rx)\tp C_c(\Rx)$.

If $f$ and $g$ are measurable functions $\Rx\to\Cplx$, let $f\cdot g$ denote their pointwise product (with the usual identifications of functions that agree~a.e.).

\begin{cor}\label{c:H-valued a.e.}
Let $\xi_1,\xi_2\in \sH$. For $h\in\Real\setminus\{0,-1\}$ let
$F(h) =\lm(1+h)\xi_1\cdot \lm(1+h^{-1})\xi_2$.
Then $F$ is equal a.e.\ to a strongly measurable, (Bochner-)square integrable function $\Rx\to\sH$, and
\[
\int_{\Rx} \norm{F(h)}_{\sH}^2 \,\frac{dh}{|h|} = \norm{\xi_1}_{\sH}^2 \norm{\xi_2}_{\sH}^2
\;.
\]
\end{cor}

\begin{proof}
We apply Lemma~\ref{l:fund unitary} with $X=\xi_1\tp\xi_2$.
As remarked in Section~\ref{ss:bochner background}, we may identify $V(X)$ with a function $\widetilde{F}\in L^2(\Rx;\sH)$, satisfying $\norm{\widetilde{F}}^2 = \norm{X}^2=\norm{\xi_1}_{\sH}^2\norm{\xi_2}_{\sH}^2$ and $\widetilde{F}(h)(t)=V(X)(h,t)$ for a.e.~$(h,t)\in \Rx\times\Rx$. The rest follows from the definition of~$V$.
\end{proof}

Note that {\it a priori}, one only expects the pointwise product of two functions in $\sH$ to lie in $L^1(\Rx)$.
The corollary shows that in fact, $F(h)\in \sH$ for a.e.\ $h\in\Rx$.
%% and that after modification on a set of measure zero, $h\mapsto F(h)$ is strongly measurable as an $\sH$-valued function.
%
In general one cannot expect $F(h)\in\sH$ for all $h\in\Rx$, since
\begin{equation}\label{eq:bad values}
\norm{F(1)}_{\sH}^2
= \int_{\Rx} \abs{ \xi_1(t/2)\xi_2(t/2)}^2 \, \frac{dt}{|t|}
= \int_{\Rx} \abs{ \xi_1(t)\xi_2(t)}^2 \, \frac{dt}{|t|}
= \norm{\xi_1\cdot\xi_2}_{\sH}^2
\end{equation}
and so taking e.g.\ $\xi_1(t)=\xi_2(t) = \mathbf{1}_{t >1} (t-1)^{-1/3}$ one sees that the RHS can be infinite.

\begin{prop}[Explicit fusion for coefficient functions of $\pi$]
\label{p:fusion of coeff}
Let $\xi_i,\eta_i\in\sH$ for $i=1,2$. Then
\[
(\xi_1 *_{\pi} \eta_1) \cdot (\xi_2*_{\pi} \eta_2) =
\int_{\Rx}  \left[ \lm(1+h)\xi_1 \cdot \lm(1+h^{-1})\xi_2 \right]
*_\pi
\left[ \lm(1+h)\eta_1 \cdot \lm(1+h^{-1})\eta_2 \right]
\,\frac{dh}{|h|}
\]
defined as the Bochner integral of an $\FA_\pi(\affR)$-valued function.
\end{prop}

\begin{proof}
Let $F(h)= \lm(1+h){\xi_1}\cdot\lm(1+h^{-1}){\xi_2}$ and
$G(h) =
\lm(1+h){\eta_1}\cdot\lm(1+h^{-1}){\eta_2}
$.
By Corollary~\ref{c:H-valued a.e.}, $F$ and $G$ are (after modification on a null subset of $\Rx$) strongly measurable as functions $\Rx\to\sH$, and square integrable (with respect to Haar measure on $\Rx$).

Therefore, the function
$h\mapsto F(h)*_\pi G(h)$ is strongly measurable and a.e.\ $\FA_\pi(\affR)$-valued; it is Bochner integrable (with respect to Haar measure on $\Rx$), since
\[ \begin{aligned}
\int_{\Rx} \norm{F(h)*_\pi G(h)}_{\FA} \,\frac{dh}{|h|}
& \leq
\int_{\Rx} \norm{F(h)}_{\sH} \norm{G(h)}_{\sH} \,\frac{dh}{|h|} \\
& \leq
\left( \int_{\Rx} \norm{F(h)}_{\sH}^2 \,\frac{dh}{|h|}\right)^{1/2}
\left( \int_{\Rx} \norm{G(h)}_{\sH}^2 \,\frac{dh}{|h|}\right)^{1/2} \\
& = \norm{\xi_1}_{\sH}\norm{\xi_2}_{\sH}\norm{\eta_1}_{\sH}\norm{\eta_2}_{\sH} \;,
\end{aligned} \]
where the final equality follows by using Corollary~\ref{c:H-valued a.e.} again.
Unpacking the definitions of $F$ and $G$, and comparing them with~\eqref{eq:fusion step1}, we see that
\[
(\xi_1 *_{\pi} \eta_1)(b,a) \,
(\xi_2 *_{\pi} \eta_2)(b,a)
=
\int_{\Rx} [F(h)*_\pi G(h)](b,a) \,\frac{dh}{|h|} \qquad\text{for all $(b,a)\in\affR$}
\]
as claimed.
\end{proof}

\begin{rem}
Our direct route to the key formula~\eqref{eq:fusion step1} relied on  {\itshape ad hoc} manipulations of integrals.
There is a more conceptual approach, based on constructing an explicit intertwining map between $\pi\tp\pi$ and $I_{\sH}\tp\pi$.
This intertwining map emerges naturally from considering the
representation $\Pi$ defined in~\eqref{eq:induced} and its description as an induced representation; details
are given in Appendix~\ref{app:induced}.
In fact, this approach was originally how we came up with the formula \eqref{eq:fusion step1}, and it motivates the technique used in Lemma~\ref{l:fund unitary}.
\end{rem}

\end{subsection}

\begin{subsection}{Defining dual convolution}\label{ss:define dual conv}
The formula in Proposition~\ref{p:fusion of coeff} immediately suggests how to define the dual convolution of two rank-one tensors in $\cS_1(\sH)=\sH\ptp\sHbar$:
given $\xi,\xi'\in\sH$ and $\eta, \eta'\in\sHbar$,
%%%% MG's version, with minor LaTeX changes by YC
\begin{equation}\label{eq:fusion of rank1}
\begin{gathered}
 (\xi\tp\eta) \boxtimes (\xi'\tp\eta') \\
 \defeq
 \int_{\Real^\times}   \left( \lambda(1+h)\xi\cdot \lambda(1+h^{-1})\xi' \right) \tp  \left( \lambda(1+h)\eta\cdot \lambda(1+h^{-1})\eta'\right) \;\frac{dh}{|h|}
 \;,
\end{gathered}
\end{equation}
%%%%
where the right-hand side is defined as a Bochner integral of a function $\Rx\to \sH\ptp\sHbar$.
The proof that this function is Bochner integrable is essentially the same as the argument used in proving Proposition~\ref{p:fusion of coeff}, so we shall not repeat it here; we record for reference that the same calculation yields the upper bound
\begin{equation}\label{eq:rank one bound}
\begin{gathered}
\int_{\Real^\times}   \norm{ \lambda(1+h)\xi\cdot \lambda(1+h^{-1})\xi' }_{\sH} \, \norm{ \lambda(1+h)\eta\cdot \lambda(1+h^{-1})\eta'}_{\sH} \,\frac{dh}{|h|}
\\
\leq
\norm{\xi}_{\sH}
\norm{\xi'}_{\sH}
\norm{\eta}_{\sH}
\norm{\eta'}_{\sH}
\;.
\end{gathered}
\end{equation}

\begin{rem}[Technical caveats]\label{r:C_c}
Strictly speaking, the integrand in \eqref{eq:fusion of rank1} is only a.e.\ $\sH\ptp\sHbar$-valued (c.f.~Equation \ref{eq:bad values}), and the null set of ``bad values'' of $h$ might depend on all four of the vectors $\xi,\xi',\eta,\eta'$.
However, one can ignore such technicalities if $\xi,\xi',\eta,\eta'\in C_c(\Rx)$. For, under this assumption, $\lm(a)\xi\cdot\xi'$  vanishes identically whenever $|a|$ is sufficiently small or sufficiently large. It follows (using continuity of translation in $\sH$ and in $C_0(\Rx)$) that the integrand in \eqref{eq:fusion of rank1} is a continuous, compactly supported  function $\Rx \setminus \{-1\}\to \sH\ptp \sHbar$, with no need to worry about various formulas holding only~a.e.
\end{rem}

We can now extend the operation $\boxtimes$ by linearity and continuity to a contractive bilinear map
$\cS_1(\sH)\times \cS_1(\sH)\to\cS_1(\sH)$,
by representing elements of $\cS_1(\sH)$ as absolutely convergent sums of rank-one tensors.
To see that this extension is well-defined and independent of how we represent elements of $\cS_1(\sH)$, note that $(\xi,\eta,\xi',\eta')\mapsto (\xi\tp\eta)\boxtimes (\xi'\tp\eta')$ defines a contractive multilinear map from $\sH\times\sHbar\times\sH\times\sHbar$ to $\sH\ptp\sHbar$,
and so by the universal property of $\ptp$, it extends uniquely to a contractive linear map
\[
\cS_1(\sH)\ptp\cS_1(\sH) = \sH\ptp\sHbar \ptp \sH\ptp\sHbar \longrightarrow \sH\ptp\sHbar =\cS_1(\sH)\;.
\]

%%%%%%%%%%%%%%%%%%%%%%%%%%%%%%%%%

\paragraph{An alternative integral formula.}
One can rewrite the defining formula \eqref{eq:fusion of rank1} as
\begin{equation}\label{eq:fusion rewrite}
(\xi\tp\eta)\boxtimes(\xi'\tp\eta')
=
\int_{\Real}( \lm(1-u)^{-1}\xi \cdot \lm(u)^{-1}\xi' )
\tp ( \lm(1-u)^{-1}\eta \cdot \lm(u)^{-1}\eta' )
\;\frac{du}{\abs{1-u} \abs{u}}
\end{equation}
after a change of variables\footnotemark\ $1-u = (1+h)^{-1}$.
\footnotetext{Changes of variables for Bochner integrals can be easily justified by verifying first for simple functions, and then passing to the limit.}
Similar comments as in Remark~\ref{r:C_c} also apply here: for instance, if $\xi,\xi',\eta,\eta'\in C_c(\Rx)$, then the integrand in \eqref{eq:fusion rewrite} is continuous from $\Real\setminus\{0,1\}$ to $\sH\ptp\sHbar$ with compact support.

Equation~\eqref{eq:fusion rewrite} should be compared with the initial calculations in \eqref{eq:getting started}.
In fact, many of the preceding results could have been formulated without the change of variables in \eqref{eq:fusion step1}.
Both formulations of dual convolution seem to be natural and useful: the formula \eqref{eq:fusion of rank1} is more closely related to the underlying general principles concerning tensor products of induced representations;
but \eqref{eq:fusion step1} is more enlightening for certain calculations, such as \eqref{eq:halfway} below.

\paragraph{An abstract definition of $\boxtimes$.}
An alternative way to think of our construction of $\boxtimes$, viewed as a bounded linear map from $\cS_1(\sH)\ptp \cS_1(\sH)$ to $\cS_1(\sH)$, is by constructing it as the composition of the maps shown in Figure~\ref{fig:fancy POV}.

\begin{figure}[htp]
\[
\begin{aligned}
(\sH \ptp \sHbar ) \ptp (\sH\ptp \sHbar)
&
\strut \xrightarrow{\rm shuffle} \strut &
(\sH \ptp \sH ) \ptp (\sHbar\ptp \sHbar)
 \\
&
\strut \xrightarrow{\rm embed} \strut  &
(\sH \tp^2 \sH ) \ptp (\sHbar\tp^2 \sHbar)
 \\
&
\strut \xrightarrow{V \tp V} \strut  &
(\sH \tp^2 \sH ) \ptp (\sHbar\tp^2 \sHbar)
 \\
&
\strut \xrightarrow{\rm identify} \strut  &
L^2(\Rx; \sH)  \ptp \overline{L^2(\Rx ; \sH)}
\\
&
\strut \xrightarrow{\rm diagonal} \strut  &
L^1(\Rx; \sH  \ptp \sHbar)
\\
&
\strut \xrightarrow{\rm trace} \strut  &
\sH  \ptp \sHbar
\end{aligned}
\]
\caption{Dual convolution as a composition of simpler operations}
\label{fig:fancy POV}
\end{figure}

We now explain briefly what each of these maps is.
\begin{itemize}
\item The ``shuffle'' map interchanges the second and third factors in the tensor product, i.e.~it sends $\xi\tp\eta\tp\xi'\tp\eta'$ to $\xi\tp\xi'\tp\eta\tp\eta'$.
\item The  ``embed'' map is self-explanatory, and $V$ is from Lemma \ref{l:fund unitary}. The map ``identify'' is the same identification described in Section~\ref{ss:bochner background} and used in Corollary~\ref{c:H-valued a.e.}.
\item
The ``diagonal'' map is given as follows: for Banach spaces $E_1$ and $E_2$ there is a canonical contraction
\[ L^2(\Rx;E_1)\ptp \overline{L^2(\Rx;E_2)} \to L^1(\Rx; E_1 \ptp \overline{E_2}) \]
which sends $F\tp G$ to $h\mapsto F(h)\tp G(h)$.
\item
The ``trace'' map is given as follows: for a Banach space $E$ there is a canonical contraction $L^1(\Rx; E)\to E$ which sends a function $F\in L^1(\Rx;E)$ to $\int_{\Rx} F$. (If we identify $L^1(\Rx;E)$ with $L^1(\Rx)\ptp E$, then the trace map is the same as slicing in the first variable against the constant function ${\bf 1}\in L^\infty(\Rx)$.)
\end{itemize}

The advantage of this approach is that all issues concerning strong measurability, or showing that various maps are well-defined and do not depend on how an element of $\sH\ptp\sHbar$ is represented as an infinite sum of tensors, are automatically taken care of by the formal identifications between various Banach spaces.
Moreover, this approach also generalizes easily to the $L^p$-setting,
or to settings with additional operator space structure.
The disadvantage is that this definition of $\boxtimes$ is rather abstract, and is less suited to concrete calculations. 

\end{subsection}
%%%%%%%%%%%%%%%%%%%%%%%%%%%%%%%%%%%%% 

\begin{subsection}{Basic properties of dual convolution}\label{ss:basic properties}
Clearly $\boxtimes$ is commutative: this follows directly from a change of variable $h\mapsto h^{-1}$ in (\ref{eq:fusion of rank1}). 
Proving that $\boxtimes$ is associative requires more work.
(Recall that even when considering the usual convolution of two $L^1$-functions on a locally compact group $G$, checking associativity directly by attempting to interchange integrals requires careful use of Fubini's theorem to justify treating identities that only hold a.e.\ as if they hold everywhere.)

To show that $\boxtimes$ is associative, it
 suffices by linearity and continuity to show that $(T_1\boxtimes T_2)\boxtimes T_3 = T_1\boxtimes (T_2\boxtimes T_3)$ when $T_i=\xi_i\tp\eta_i$ for $\xi_i,\eta_i \in C_c(\Rx)$ ($i=1,2,3$).
In the following calculations, we shall adopt the following notational convention to make our formulas more manageable. Given a function in $C_c(\Rx)$ which is obtained from $\xi_1$, $\xi_2$ and $\xi_3$ by some explicit formula $\cT[\xi_1,\xi_2,\xi_3]$, we shall write
\REPEAT{\cT[\xi_1,\xi_2,\xi_3]}
to mean
\[ \cT[\xi_1,\xi_2,\xi_3] \tp \cT[\eta_1,\eta_2,\eta_3] \in
{C_c(\Rx)\tp C_c(\Rx).} \]

Since $(\xi_1\tp\eta_1)\boxtimes (\underline{\quad})$ is a bounded linear map $\sH\ptp\sHbar \to \sH\ptp\sHbar$, it commutes with the Bochner integral. In particular,
%For any bounded linear map $S: \sH\ptp\sHbar \to \sH\ptp\sHbar$
\begin{equation}\label{eq:iterated}
\begin{aligned}
& \phantom{\quad}
(\xi_1\tp\eta_1)\boxtimes \Bigl(
(\xi_2\tp\eta_2)
\boxtimes
(\xi_3\tp\eta_3)
\Bigr)
\\
& =
\int_{\Real}
(\xi_1\tp\eta_1)\boxtimes \Bigl(
\ \REPEAT{ \lambda(1-u)^{-1}\xi_2\cdot \lambda(u)^{-1}\xi_3 } \
\Bigr)
\;\frac{du}{\abs{1-u}\abs{u}}
\\
& =
\int_{\Real}
\left(
\int_{\Real}
\REPEAT{
\lambda(1-v)^{-1}\xi_1 \cdot \lambda(v)^{-1}\left[ \lambda(1-u)^{-1}\xi_2\cdot \lambda(u)^{-1}\xi_3 \right]}
\;\frac{dv}{\abs{1-v}\abs{v}}
\right)
\frac{du}{\abs{1-u}\abs{u}}
\\
& =
\int_{\Real}
\left(
\int_{\Real}
\REPEAT{
\lambda(1-v)^{-1}\xi_1 \cdot \lambda(v-uv)^{-1}\xi_2\cdot \lambda(uv)^{-1}\xi_3 }
\;\frac{dv}{\abs{1-v}\abs{v}}
\right)
\frac{du}{\abs{1-u}\abs{u}}
\\
\end{aligned}
\end{equation}

The expression in the inner integral is measurable and Bochner integrable as a function $\Real^2\to \sH\ptp\sHbar$ (since it is continuous with compact support and vanishes in a neighbourhood of $\{0,1\}\times\{0,1\}$). So by Fubini's theorem for Bochner integrals
(see e.g. \cite[Theorem B.41]{williams}),
 we may rewrite \eqref{eq:iterated} as a double integral and perform succesive changes of variables $u \mapsto u/v$, $v\mapsto 1-v$ to obtain
\begin{equation}\label{eq:halfway}
\begin{gathered}
(\xi_1\tp\eta_1)\boxtimes \Bigl(
(\xi_2\tp\eta_2)
\boxtimes
(\xi_3\tp\eta_3)
\Bigr)
\\
=
\int_{\Real^2}
\REPEAT{
\lambda(1-v)^{-1}\xi_1 \cdot \lambda(v-u)^{-1}\xi_2\cdot \lambda(u)^{-1}\xi_3 }
\;\frac{d(v,u)}{\abs{1-v}\abs{v-u}\abs{u}}
\\
=
\int_{\Real^2}
\REPEAT{
\lambda(v)^{-1}\xi_1 \cdot \lambda(1-v-u)^{-1}\xi_2\cdot \lambda(u)^{-1}\xi_3 }
\;\frac{d(v,u)}{\abs{v}\abs{1-v-u}\abs{u}}
\\
\end{gathered}
\end{equation}
One can use similar arguments to expand
$\Bigl(
(\xi_1\tp\eta_1)
\boxtimes
(\xi_2\tp\eta_2)
\Bigr)
\boxtimes (\xi_3\tp\eta_3)$ as a double integral with values in $\sH\ptp\sHbar$, and show by appropriate changes of variable that this is equal to the right-hand side of \eqref{eq:halfway}. Alternatively, observe that since $\boxtimes$ is commutative,
\[
\Bigl(
(\xi_1\tp\eta_1)
\boxtimes
(\xi_2\tp\eta_2)
\Bigr)
\boxtimes (\xi_3\tp\eta_3) 
=
(\xi_3\tp\eta_3)
\boxtimes
\Bigl(
(\xi_2\tp\eta_2)
\boxtimes
(\xi_1\tp\eta_1)
\Bigr);
\]
then observe that the value of the last integral in Equation \eqref{eq:halfway} is unchanged if one swaps $\xi_1\tp\eta_1$ with $\xi_3\tp\eta_3$ (since this corresponds to interchanging the variables $u$ and $v$ in the integral).

Note that Proposition~\ref{p:fusion of coeff} can be rephrased as
\begin{equation}
\Psi((\xi\tp\eta)\boxtimes(\xi'\tp\eta')) = \Psi(\xi\tp\eta)\Psi(\xi'\tp\eta'),
\end{equation}
and so by linearity and continuity, it follows that $\Psi(T\boxtimes T')=\Psi(T)\Psi(T')$ for all $T,T'\in\cS_1(\sH)$.
This gives an independent proof that $\FA_\pi(\affR)$ is closed under pointwise product. 
It could also have been used to prove associativity of $\boxtimes$, by transferring it from associativity of pointwise product in $\FA_\pi(\affR)$. We believe that the direct proof given above has independent interest, especially in light of the symmetry displayed by the formula in Equation~\eqref{eq:halfway}.

\paragraph{Summary.}
We sum up the results of this section in the following theorem.

\begin{thm}[Dual convolution on $\cS_1(\sH)$]
The operation $\boxtimes$, defined on pairs of elementary tensors
by the formula \eqref{eq:fusion of rank1},
extends to a contractive bilinear map $\cS_1(\sH)\times \cS_1(\sH)\to\cS_1(\sH)$, which makes $\cS_1(\sH)$ into a commutative Banach algebra. If we denote this Banach algebra by $\cA$, then $\Psi:\cA\to\FA_\pi(\affR)$ is an isometric isomorphism of Banach algebras.
\end{thm}

\end{subsection}

\begin{subsection}{Dual convolution at the level of functions}\label{ss:view as functions}
Trace-class operators on $\sH=L^2(\Rx)$ are often given not as explicit sums of rank-one tensors, but as integral operators defined by certain kernel functions $\Rx\times\Rx\to\Cplx$. In this section we provide a description of dual convolution that may be easier to apply in such cases.

We may view elements of $\cS_1(\sH)$ as measurable functions on $\Rx\times\Rx$,
 as follows.
First note that complex conjugation of functions defines a $\Cplx$-linear isometric isomorphism of vector spaces from $\overline{L^2(\Rx)}$ onto $L^2(\Rx)$, which extends to an isometric isomorphism
\begin{equation}
\tilde{\iota}: L^2(\Rx)\ptp \overline{L^2(\Rx)} \to L^2(\Rx)\ptp L^2(\Rx)
\qquad; \qquad  \xi\tp\eta \mapsto \xi\tp\overline{\eta}.
\end{equation}
Furthermore, the natural map $L^2(\Rx)\ptp L^2(\Rx) \to L^2(\Rx)\tp^2 L^2(\Rx)$
is linear and norm-decreasing, and it is injective since Hilbert spaces have the approximation property.
Finally, note that we may identify $L^2(\Rx)\tp^2 L^2(\Rx)$ with $L^2(\Rx\times \Rx)$.

Thus, up to a.e.\ equivalence\footnotemark\ we can view any $T\in\cS_1(\sH)$ as a measurable function on $\Rx\times\Rx$ which is square-integrable (with respect to the measure $|s|^{-1}|t|^{-1}\,d(s,t)$).
\footnotetext{There is a subtler notion available when viewing elements of $\cS_1(L^2(\Omega))$ as functions on $\Omega\times\Omega$; rather than quotienting out by the equivalence relation ``agree except on a null subset of $\Omega\times\Omega$'', one uses the finer equivalence relation ``agree except on a \emph{marginally null} subset''. This notion, which orginates in pioneering work of Arveson on operator synthesis, is not needed for our paper.}
%%
%% This function is the integral kernel corresponding to the operator~$T$;
For ease of notation, we shall denote this function also by~$T$, suppressing mention of the embedding~$\tilde{\iota}$.
With this convention,
\[
\ip{T\alpha}{\beta} = \int_{\Rx\times \Rx}  \overline{\beta(s)} T(s,t) \alpha(t)
\; \frac{d(s,t)}{\abs{s}\abs{t}} \qquad(\alpha,\beta\in\sH),
\]
which is the usual form in which an integral operator is given.
{\bf Warning:} with this convention, if $T$ is a rank-one tensor $\xi\tp\eta \in \sH\ptp\sHbar$ then $T(s,t)=\xi(s)\overline{\eta(t)}$ for $s,t\in\Rx$.

\begin{prop}[Pointwise formulas for dual convolution]
\label{p:pointwise dc}
Let $T_1$, $T_2\in\cS_1(\sH)=\sH\ptp\sHbar$. Then for a.e.~$(s,t)\in\Rx\times\Rx$
\begin{subequations}
\begin{align}
(T_1\boxtimes T_2)(s,t)
& = \label{eq:fletcher1}
    \int_{-\infty}^\infty
T_1\left( \frac{s}{1+h}, \frac{t}{1+h}\right)
T_2\left( \frac{hs}{1+h}, \frac{ht}{1+h}\right)
\,\frac{dh}{|h|} \\
& = \label{eq:fletcher2}
    \int_{-\infty}^\infty
T_1((1-u)s, (1-u)t )
T_2(us,ut)
\,\frac{du}{\abs{1-u}\abs{u}}
\end{align}
\end{subequations}
where both of the integrals above are absolutely convergent for a.e.~$(s,t)\in\Rx\times\Rx$.
\end{prop}
\begin{proof}
When $T_1$ and $T_2$ are rank-one tensors, this follows from the definition of $\boxtimes$.
Hence it is true when $T_1$ and $T_2$ are finite rank operators. 
Every trace class operator is the limit in trace-norm of finite rank operators, and by going down to a subsequence we can assume that the convergence holds pointwise a.e.

Now observe that if $T=\sum_{n=1}^\infty f_n\otimes g_n$ where $\sum_{n=1}^\infty \norm{f_n}_{\sH} \norm{g_n}_{\sH}<\infty$, the trace-class operator $R=\sum_{n=1}^\infty |f_n|\otimes |g_n|$ satisfies $R(s,t)\geq |T(s,t)|$~a.e.
The result now follows using the Lebesgue dominated convergence theorem,
replacing $T_1$ and $T_2$ in \eqref{eq:fletcher1} or \eqref{eq:fletcher2} with ``dominating operators'' $R_1$ and~$R_2$.
\end{proof}

%%%%%%%%%%%%%%%%%%%%%%%%%%%%

\end{subsection}
\end{section}

\begin{section}{An explicit derivation from $\cA$ to its dual}
\label{s:derivation}

In this section we construct an explicit derivation $D:\cA\to\cA^*$ and study some of its operator-theoretic properties.
We will relate $D$ to constructions in \cite{Choi-Ghandehari14} in the next section.

We briefly review some general definitions. For a Banach algebra $A$, each $\Phi\in (A\ptp A)^*$ corresponds to a bounded linear map  $A\to A^*$ defined by $a\mapsto \Phi(a \tp\overline{\quad})$.
This map $A\to A^*$ is a {\itshape derivation} if the following identity holds:
\begin{equation}\label{eq:leibniz}
\Phi(a_1a_2\tp a_0) = \Phi(a_2\tp a_0a_1) + \Phi(a_1\tp a_2a_0)
\qquad\text{for all $a_0,a_1,a_2\in A$.}
\end{equation}
The derivation is said to be {\itshape cyclic} if $\Phi(a\tp b)=-\Phi(b\tp a)$ for all $a,b\in A$.

\begin{dfn}
Given $\xi\in\sH$, let $S\xi(t)=\sgn(t)\xi(t)$ and $R\xi(t)=\xi(-t)$,
where $\sgn$ is the sign function $\Rx\to \{\pm 1\}$, $t\mapsto \frac{t}{\abs{t}}$.
Clearly $S$ and $R$ are isometric, linear involutions on~$\sH$.
\end{dfn}
%
%Clearly, $K^{-1}L$ extends to the bounded operator $S$ on $\sH$.

Although we do not consider coefficient functions in this section, note that for every $\xi,\eta\in\sH$ we have $\overline{\xi*_\pi\eta}=R\overline{\xi}*_\pi R\overline{\eta}$.

\paragraph{Constructing our derivation.}
Define a multilinear map $\Phi:\sH\ptp\sHbar \ptp \sH\ptp\sHbar \to\Cplx$ by 
\begin{equation}\label{eq:define derivation}
\Phi(\xi_1\otimes \eta_1\ \otimes \xi_0\otimes \eta_0)
\defeq
\ip{S\xi_1}{R\overline{\xi_0}} \ip{R\overline{\eta_0}}{\eta_1}.
\end{equation}
We view $\Phi$ as a bilinear form on $\cA$, and define $D$ to be the corresponding operator $\cA\to\cA^*$.

For $i=0,1$, if we write $T_i=\xi_i\tp\eta_i$ and use the convention $T_i(s,t)=\xi_i(s)\overline{\eta_i(t)}$ as in Section~\ref{ss:view as functions}, then we can rewrite \eqref{eq:define derivation}~as:
\begin{equation}\label{eq:der easier formula}
\Phi(T_1\tp T_0) = \int_{\Rx\times\Rx} \sgn(s) T_1(s,t) T_0(-s,-t) \, \frac{d(s,t)}{|s||t|}.
\end{equation}

\begin{prop}\label{p:D-is-cyclic-and-wc}
$D:\cA\to\cA^*$ is cyclic and weakly compact.
\end{prop}

\begin{proof}
The identity \eqref{eq:der easier formula} shows that $\Phi(T_1\tp T_0)=-\Phi(T_0\tp T_1)$ when $T_0$ and $T_1$ are rank-one tensors; by linearity and continuity it holds for all $T_0,T_1\in\cA$.

It also follows from \eqref{eq:der easier formula}, using the Cauchy--Schwarz inequality, that $\Phi$ extends to a bounded bilinear form on $L^2(\Rx\times\Rx)$. Hence the operator $D:\cA\to\cA^*$ factors through the embedding of $\cA$ into $L^2(\Rx\times\Rx)$; in particular $D$ is weakly compact.
\end{proof}

Next, we show that $D$ is a derivation by showing that $\Phi$ satisfies the identity~\eqref{eq:leibniz}.
 
\begin{thm}[Derivation identity]
For every $T_1, T_2, T_0\in \cA$, we have 
\begin{equation}\label{eq:der identity}
\Phi((T_1\boxtimes T_2)\otimes T_0)=\Phi(T_2\otimes (T_0\boxtimes T_1))+\Phi(T_1\otimes (T_2\boxtimes T_0)).
\end{equation}
\end{thm}

\begin{proof}
% In what follows, to simplify notation, we omit the domain of integration: it will be either $\Rx$, $(\Rx)^2$ or $(\Rx)^3$ depending on context.
%
By linearity and continuity, it suffices\footnotemark\ to verify \eqref{eq:der identity} in the special case where $T_0$, $T_1$ and $T_2$ are rank-one tensors in $C_c(\Rx)\tp C_c(\Rx)$.
\footnotetext{Strictly speaking, this reduction step is not necessary, but it removes any need to consider technicalities about interchanging the order of various integrals that now follow.}

We now consider the three terms in \eqref{eq:der identity}, using \eqref{eq:der easier formula} and Fubini's theorem. In each case the integral is taken over $(\Rx)^3$:

\[ \begin{aligned}
& \Phi((T_1\boxtimes T_2)\tp T_0)
\\
& = \int \sgn(s)
\, T_1\bigl(\frac{s}{1+h}, \frac{t}{1+h}\bigr)
\, T_2\bigl(\frac{hs}{1+h}, \frac{ht}{1+h}\bigr)
\, T_0(-s,-t)\; \frac{d(h,s,t)}{|h||s||t|}
\;.
\end{aligned}
\]
Also,
\[
\begin{aligned}
& \Phi(T_2\tp (T_0\boxtimes T_1))
\\
& =
\int \sgn(s)
\, T_1\bigl(\frac{-hs}{1+h},\frac{-ht}{1+h}\bigr)
\, T_2\bigl(s,t)
\, T_0\bigl(\frac{-s}{1+h},\frac{-t}{1+h}\bigr) 
\frac{d(h,s,t)}{|h||s||t|}
\\
& = \int \sgn((1+h)s)
\, T_1\bigl(-hs,-ht\bigr)
\, T_2\bigl((1+h)s, (1+h)t\bigr)
\, T_0(-s,-t)\; \frac{d(h,s,t)}{|h||s||t|}
\;,
\end{aligned}
\]
where the last equality used the change of variables $s\mapsto (1+h)s$, $t\mapsto (1+h)t$.
A similar calculation yields
\[
\begin{aligned}
& \Phi(T_1\tp(T_2\boxtimes T_0))
\\
& = \int \sgn( (1+\tfrac{1}{h})s) \, T_1\bigl((1+\tfrac{1}{h})s,(1+\tfrac{1}{h})t\bigr)
\, T_2\bigl(-\tfrac{1}{h}s, -\tfrac{1}{h}t\bigr)
\, T_0(-s,-t)\; \frac{d(h,s,t)}{|h||s||t|}
\;.
\end{aligned}
\]

For every $s,t\in\Rx$, define
\begin{eqnarray*}
  \ONE(s,t)&\defeq &
  \int_{\Real^\times} \sgn(s)\, T_1\bigl(\frac{s}{1+h} , \frac{t}{1+h}\bigr) T_2 \bigl( \frac{hs}{1+h} , \frac{ht}{1+h}\bigr) \;\frac{dh}{|h|}\\
  \TWO(s,t)&\defeq &
  \int_{\Real^\times}\sgn((1+h)s) T_1( -hs, -ht)  T_2\bigl( (1+h)s, (1+h)t)\;\frac{dh}{|h|}\\
  \THREE(s,t)&\defeq &
  \int_{\Real^\times} \sgn( (1+\tfrac{1}{h}) s)\, T_1\bigl( (1+\tfrac{1}{h})s ,(1+\tfrac{1}{h})t \bigr) T_2\bigl(-\tfrac{1}{h}s, -\tfrac{1}{h}t \bigr) \;\frac{dh}{|h|}
\end{eqnarray*}
To prove that \eqref{eq:der identity} holds, it suffices to show that $\ONE(s,t)= \TWO(s,t)+ \THREE(s,t)$  for (almost) every $s,t\in\Real^\times$.
For $\TWO(s,t)$: the change of variables $h\mapsto -\frac{1}{1+h}$ sends $\frac{dh}{|h|}$ to $\frac{dh}{\abs{1+h}}$ and sends $1+h$ to $1-\frac{1}{1+h}=\frac{h}{1+h}$, so that
\begin{equation}
\label{eq:II in place}
\begin{aligned}
\TWO(s,t)
&= \int_{\Rx} \sgn\left(\frac{hs}{1+h}\right) \frac{|h|}{|1+h|} \,
T_1\bigl(\frac{s}{1+h}, \frac{t}{1+h}\bigr)
\,T_2 \bigl( \frac{hs}{1+h} , \frac{ht}{1+h}\bigr)\;\frac{dh}{|h|} \\
&= \int_{\Rx} \frac{1}{|s|} \frac{hs}{1+h}  \,
T_1\bigl(\frac{s}{1+h}, \frac{t}{1+h}\bigr)
\,T_2 \bigl( \frac{hs}{1+h} , \frac{ht}{1+h}\bigr)\;\frac{dh}{|h|} \\
\end{aligned}
\tag{$**$}
\end{equation}

For $\THREE(s,t)$: the change of variables $h \mapsto -(1+\tfrac{1}{h})$ sends $\frac{dh}{|h|}$ to $\frac{dh}{\abs{1+h}\abs{h}}$
and sends $1+\tfrac{1}{h}$ to $1 - \frac{1}{1+h^{-1}}= \frac{1}{1+h}$,
so that
\begin{equation}
\label{eq:III in place}
\begin{aligned}
\THREE(s,t)
& = \int_{\Rx} \sgn\left(\frac{s}{1+h}\right) \frac{1}{|1+h|} \,
T_1\bigl(\frac{s}{1+h}, \frac{t}{1+h}\bigr)
\, T_2 \bigl( \frac{hs}{1+h} , \frac{ht}{1+h}\bigr)\;\frac{dh}{|h|}
\\
&= \int_{\Rx} \frac{1}{|s|} \frac{s}{1+h}  \,
T_1\bigl(\frac{s}{1+h}, \frac{t}{1+h}\bigr)
\, T_2 \bigl( \frac{hs}{1+h} , \frac{ht}{1+h}\bigr)\;\frac{dh}{|h|} \\
\end{aligned}
\tag{$***$}
\end{equation}
Adding \eqref{eq:II in place} and \eqref{eq:III in place} and recalling that $\frac{s}{|s|}=\sgn(s)$, we obtain $\ONE(s,t)$ as required.
\end{proof}

Finally, we show that $D:\cA\to\cA^*$ is completely bounded after composing with the transpose map on $\cA^*=\Bdd(\sH)$.

\begin{dfn}[Transpose operator on $\Bdd(\sH)$]\label{d:transpose}
Define a \emph{linear} isometry $\top_*:\sH\ptp\sHbar\to \sH\ptp\sHbar$ by $\top_*(\xi\tp\eta)= \overline{\eta}\tp\overline{\xi}$,
and let $\top=(\top_*)^* : (\sH\ptp\sHbar)^*\to (\sH\ptp\sHbar)^*$.
We call $\top$ the {\itshape transpose operator}, since if we identify $(\sH\ptp\sHbar)^*$ with $\Bdd(\sH)$ we have
\[
\overline{\top(b)\xi}=b^*(\overline{\xi}) \qquad(b\in\Bdd(\sH), \xi\in\sH).
\]
\end{dfn}

To verify complete boundedness of $\top D:\cA\to\cA^*$ we use the following characterization.
Recall that we have natural injective maps
\[
(\sH \ptp \sHbar ) \ptp (\sH\ptp \sHbar)
 \xrightarrow{\rm shuffle} 
(\sH \ptp \sH ) \ptp (\sHbar\ptp \sHbar)
\xrightarrow{\rm embed}
(\sH \tp^2 \sH ) \ptp (\sHbar\tp^2 \sHbar)
\]
where ``shuffle'' swaps the second and third factors in the tensor product. Now define
\[
s = {\rm embed}\circ{\rm shuffle}: \cS_1(\sH) \ptp \cS_1(\sH) \to \cS_1(\sH\tp^2\sH).
\]
Then, for a given $\Lambda \in ( \cS_1(\sH)\ptp\cS_1(\sH))^*$, the corresponding map $\cS_1(\sH)\to \cS_1(\sH)^*$ is {\itshape completely bounded} (with respect to the natural operator space structure on $\cS_1(\sH)$) if and only if $\Lambda \circ s^{-1}$ extends continuously to a bounded linear functional on $\cS_1(\sH\tp^2\sH)$.
(This may be taken as a
working definition of complete boundedness in this special case; it can also be derived from
facts about row and column Hilbert spaces and operator space tensor products, see e.g.\ \cite[Corollary 7.1.5 and Proposition 7.2.1]{ER_OSbook}.)

\begin{proof}[Proof that $\top D$ is completely bounded]
By \eqref{eq:define derivation}, $\top D: \cA\to \cA^*$ corresponds to the linear functional $\Psi:\cA\ptp\cA\to\Cplx$ defined by
\begin{equation}
\Psi(\xi\tp\eta \tp \xi'\tp \eta')
=\Phi(\xi\tp\eta \tp \overline{\eta'}\tp\overline{\xi'})
= \ip{S\xi}{R\eta'} \ip{R\xi'}{\eta}\,.
\end{equation}
We have
\begin{equation}
(\Psi\circ s^{-1})( \xi\tp\xi' \tp \eta\tp\eta')
= \ip{(S\tp R)(\xi \tp \xi')}{(R\tp {\rm id})(\eta'\tp \eta)} \,;
\end{equation}
since $S\tp R$ and $R\tp {\rm id}$ are unitary operators on $\sH\tp^2\sH$, it follows that $\Psi\circ s^{-1}$ extends continuously to an element of $(\cS_1(\sH\tp^2\sH))^*$.
\end{proof}

\begin{rem}
In the language of \cite{choi_Z2-pre}, $D:\cA\to\cA^*$ is {\itshape co-completely bounded}, since $\top$ ``reverses the operator space structure'' on $\Bdd(\sH)$.
\end{rem}

\end{section}

%%%%%%%%%%%%%%%%%%%%%%%%%%%%%%%%%%%%%%%%%%%%%%%%%%%%%

\begin{section}{Dual convolution for $\affRpos$}\label{s:subalgebra}
The semidirect product $\affRpos$ may be viewed as an open subgroup of $\affR$,
by identifying it with $\{(b,a)\colon b\in\Real,a\in\Rxpos\}$.
General results on Fourier algebras of open subgroups then allow us to identify $\FA(\affRpos)$ with a closed subalgebra of $\FA(\affR)$.
More precisely, let us introduce the non-standard notation:
\begin{equation}
\FA_e(\affR) \defeq\{ f\in \FA(\affR) \colon \text{$f$ vanishes outside $\affRpos$} \}.
\end{equation}
Given $f\in\FA(\affR)$, define $P_e(f)(b,a)=f(b,a)$ for $a>0$ and $P_e(f)(b,a)=0$ for $a<0$.
Then $P_e$ is a (completely) contractive projection from $\FA(\affR)$ onto $\FA_e(\affR)$,
and the composition of the maps
\[
\FA_e(\affR) \xrightarrow{{\rm inc.}} \FA(\affR)\xrightarrow{{\rm restr.}}\FA(\affRpos)
\]
is a completely isometric bijection.
(See e.g.~\cite[Prop.~2.4.1]{KL_AGBGbook} for a summary of the necessary facts about Fourier algebras of open subgroups.)

%%Since $\Psi:\cA\to\FA(\affR)$ is a bijective homomorphism, it is natural to ask for a concrete description of the subalgebra of $\cA$ that corresponds to $\FA_e(\affR)$. This is the purpose of the current section.

We now proceed to identify the subalgebra of $\cA$ that corresponds to $\FA_e(\affR)$ and hence models $\FA(\affRpos)$.
Recall the operator $S:\sH\to\sH$ given by $(S\xi)(t)=\sgn(t)\xi(t)$.
$S$~is an isometric involution, so it has eigenvalues $\pm 1$, and $\sH$ decomposes as an orthogonal direct sum of the corresponding eigenspaces,
$\sH = \sH^+ \oplus_2 \sH^-$.
The spaces $\sH^{\pm}$ have the following explicit description:
\begin{equation}
\sH^{+} = \{ \xi \in \sH \colon \xi(t)=0 \;\text{for a.e. $t < 0$} \}
\quad,\quad
\sH^{-} = \{ \xi \in \sH \colon \xi(t)=0 \;\text{for a.e. $t > 0$} \}.
\end{equation}

Let $P^{\pm}$ be the orthogonal projection of $\sH$ onto $\sH^{\pm}$, and define $P_{\rm diag} :\sH\ptp\sHbar \to \sH\ptp\sHbar$ by
\begin{equation}
P_{\rm diag}(\xi\tp\eta) = P^{+}\xi \tp P^{+}\eta + P^{-}\xi \tp P^{-}\eta \;.
\end{equation}

\begin{lem}
$P_{\rm diag}$ is a norm-one projection.
\end{lem}

\begin{proof}
Since $P^{+}$ and $P^{-}$ are complementary projections, a direct calculation yields $(P_{\rm diag})^2 = P_{\rm diag}$. Therefore, it suffices to show that $P_{\rm diag}$ is contractive:
\[ \begin{aligned}
\norm{
P_{\rm diag}(\xi\tp\eta)}_{\sH\ptp\sHbar}
& \leq
\norm{P^{+}\xi}_{\sH} \norm{P^{+}\eta}_{\sH}
+ \norm{P^{-}\xi}_{\sH} \norm{P^{-}\eta}_{\sH}
\\
& \leq
\left( \norm{P^{+}\xi}_{\sH}^2  + \norm{P^{-}\xi}_{\sH}^2\right)^{1/2}
\left( \norm{P^{+}\eta}_{\sH}^2  + \norm{P^{-}\eta}_{\sH}^2\right)^{1/2}
  & \expo{Cauchy--Schwarz in $\Real^2$} \\
& = \norm{\xi}_{\sH} \norm{\eta}_{\sH} \;.
  & \expo{Pythagoras in $\sH$}
\end{aligned}
\]
By the definition of the projective tensor norm, it follows that $\norm{P_{\rm diag}(T)}\leq \norm{T}$ for all $T\in \sH\ptp\sHbar$.
\end{proof}

\begin{rem}
The decomposition $\sH=\sH^+\oplus_2\sH^-$ gives a decomposition of $T\in\cS_1(\sH)$ as a $2\times 2$ block matrix
\[
T = \begin{pmatrix} P_{+}TP_{+} & P_{+}TP_{-} \\ P_{-}TP_{+} & P_{-}TP_{-} 
\end{pmatrix}
\]
If we identify $\sH\ptp\sHbar$ with $\cS_1(\sH)$, then $P_{\rm diag}$ corresponds to ``compression to the diagonal''.
\end{rem}

\begin{prop}\label{prop:Adiag-direct}
 $P_{\rm diag}$ is a homomorphism.
\end{prop}
\begin{proof}
It suffices to prove that 
$P_{\rm diag}(T_1)\boxtimes P_{\rm diag}(T_2)= P_{\rm diag}(T_1\boxtimes T_2)$ for $T_1=\xi_1\tp\eta_1$ and $T_2=\xi_2\tp\eta_2$.
To simplify our formulas slightly, we write $\xi_1^+=P^{+}\xi_1$, etc., and for functions $f$, $g$ defined on the set $\{+,-\}$ we write $\sum_{\pm} f(\pm)g(\mp)$ for $f(+)g(-)+f(-)g(+)$.

We have
\begin{equation}\label{eq:four pieces}
P_{\rm diag}(T_1)\boxtimes P_{\rm diag}(T_2)
=
\left\{\begin{aligned}
\phantom{+} (\xi_1^+\tp\eta_1^+)\boxtimes(\xi_2^+\tp\eta_2^+) 
& +(\xi_1^-\tp\eta_1^-)\boxtimes(\xi_2^-\tp\eta_2^-) \\
 +(\xi_1^+\tp\eta_1^+)\boxtimes(\xi_2^-\tp\eta_2^-) 
& +(\xi_1^-\tp\eta_1^-)\boxtimes(\xi_2^+\tp\eta_2^+)
\end{aligned}
\right.
\end{equation}

To analyze each of these four terms, we consider the effect of $\lm(1+h)$ and $\lm(1+h^{-1})$ on vectors in $\sH^+$ or $\sH^{-}$, as $h$ varies over $\Rx$.
Note that if $\al\in \sH^\pm$ and $\beta\in \sH^{\mp}$ (i.e.\ $\al$ and $\beta$ have ``different parity'') then $\al\cdot\beta =0$ as elements of $\sH$.
Also: if $a>0$ then $\lm(a)(\sH^{\pm})=\sH^{\pm}$;
and if $a<0$ then $\lm(a)(\sH^{\pm})=\sH^{\mp}$\/.
Using these facts,
\begin{itemize}

\item if $h<0$ and $h\neq -1$, then precisely one of $1+h$ or $1+h^{-1}$ is negative, and so
\[ \lm(1+h)\xi_1^{\pm} \cdot \lm(1+h^{-1})\xi_2^{\pm}
= \lm(1+h)\eta_1^{\pm} \cdot \lm(1+h^{-1})\eta_2^{\pm} = 0 \;; \] 

\item if $h>0$, then $1+h$ and $1+h^{-1}$ are both positive, and so
\[ \lm(1+h)\xi_1^{\pm} \cdot \lm(1+h^{-1})\xi_2^{\mp}
= \lm(1+h)\eta_1^{\pm} \cdot \lm(1+h^{-1})\eta_2^{\mp} = 0 \;. \] 

\end{itemize}
Therefore, considering the four terms in \eqref{eq:four pieces}, we obtain
\begin{subequations}
\begin{equation}\label{eq:pm-pm}
\begin{gathered}
(\xi_1^{+}\tp\eta_1^{+})\boxtimes(\xi_2^{+}\tp\eta_2^{+}) 
+
(\xi_1^{-}\tp\eta_1^{-})\boxtimes(\xi_2^{-}\tp\eta_2^{-}) 
\\
= \sum\nolimits_{\pm}
\int_0^\infty 
\REPEAT{
\left(\lm(1+h)\xi_1^{\pm} \cdot \lm(1+h^{-1})\xi_2^{\pm}\right)
}
%\tp\left(\lm(1+h)\eta_1^{\pm} \cdot \lm(1+h^{-1})\eta_2^{\pm}\right) 
\;\frac{dh}{|h|}
\end{gathered}
\end{equation}
and
\begin{equation}\label{eq:pm-mp}
\begin{gathered}
(\xi_1^{+}\tp\eta_1^{+})\boxtimes(\xi_2^{-}\tp\eta_2^{-}) 
+
(\xi_1^{-}\tp\eta_1^{-})\boxtimes(\xi_2^{+}\tp\eta_2^{+}) 
\\
= \sum\nolimits_{\pm}
\int_{-\infty}^0 
\REPEAT{
\left(\lm(1+h)\xi_1^{\pm} \cdot \lm(1+h^{-1})\xi_2^{\mp}\right)
}
%\tp\left( \lm(1+h)\eta_1^{\pm} \cdot \lm(1+h^{-1})\eta_2^{\mp} \right)
\;\frac{dh}{|h|}
\end{gathered}
\end{equation}
\end{subequations}
where we have used the same notational convention as in Section~\ref{ss:basic properties} to simplify the formulas.

Now we consider
$P_{\rm diag}(T_1\boxtimes T_2)$; since the Bochner integral commutes with bounded linear maps, this equals
\begin{equation}\label{eq:watching chess}
\begin{gathered}
\int_{\Rx} P_{\rm diag}\left[
\REPEAT{
(\lambda(1+h)\xi_1\cdot \lambda(1+h^{-1})\xi_2)
}
%\tp (\lambda(1+h)\eta_1\cdot \lambda(1+h^{-1})\eta_2)
\right] \;
\frac{dh}{|h|} 
\end{gathered}
\end{equation}
For $-\infty < h<-1$, we have $1+h <0 < 1+h^{-1}$; hence
\[
P^{\pm}[ \lm(1+h)\xi_1\cdot \lm(1+h^{-1})\xi_2]
= \lm(1+h)\xi_1^{\mp}\cdot \lm(1+h^{-1})\xi_2^{\pm}
\]
For $-1 < h < 0$, we have $1+h > 0 > 1+h^{-1}$; hence
\[
P^{\pm}[ \lm(1+h)\xi_1\cdot \lm(1+h^{-1})\xi_2]
= \lm(1+h)\xi_1^{\pm}\cdot \lm(1+h^{-1})\xi_2^{\mp}
\]
For $0<h <\infty$, we have $1+h > 0$ and $1+h^{-1} > 0$; hence
\[
P^{\pm}[ \lm(1+h)\xi_1\cdot \lm(1+h^{-1})\xi_2]
= \lm(1+h)\xi_1^{\pm}\cdot \lm(1+h^{-1})\xi_2^{\pm}
\]
Therefore, splitting the integral in \eqref{eq:watching chess} into three pieces, and recalling that $P_{\rm diag}=\sum_{\pm} P^\pm\tp P^\pm$,
we obtain
\[
\begin{aligned}
P_{\rm diag}(T_1\boxtimes T_2)
& =\left\{
\begin{aligned}
\int_{-\infty}^{-1}
\sum\nolimits_{\pm}
\REPEAT{
 \left(\lm(1+h)\xi_1^{\mp}\cdot \lm(1+h^{-1})\xi_2^{\pm}\right)
}
% \tp\left(\lm(1+h)\eta_1^{\mp}\cdot \lm(1+h^{-1})\eta_2^{\pm}\right)
\;\frac{dh}{|h|}
\\
+\int_{-1}^{0}
\sum\nolimits_{\pm}
\REPEAT{
 \left(\lm(1+h)\xi_1^{\pm}\cdot \lm(1+h^{-1})\xi_2^{\mp}\right)
}
% \tp\left(\lm(1+h)\eta_1^{\pm}\cdot \lm(1+h^{-1})\eta_2^{\mp}\right)
\;\frac{dh}{|h|}
\\
+ \int_{0}^{\infty}
\sum\nolimits_{\pm}
\REPEAT{
 \left(\lm(1+h)\xi_1^{\pm}\cdot \lm(1+h^{-1})\xi_2^{\pm}\right)
}
% \tp\left( \lm(1+h)\eta_1^{\pm}\cdot \lm(1+h^{-1})\eta_2^{\pm}\right)
\;\frac{dh}{|h|}
\end{aligned}
\right. \\
& = \left\{
\begin{aligned}
\int_{-\infty}^{0}
\sum\nolimits_{\pm}
\REPEAT{
 \left(\lm(1+h)\xi_1^{\mp}\cdot \lm(1+h^{-1})\xi_2^{\pm}\right)
}
% \tp\left( \lm(1+h)\eta_1^{\mp}\cdot \lm(1+h^{-1})\eta_2^{\pm}\right)
\;\frac{dh}{|h|}
\\
+ \int_{0}^{\infty}
\sum\nolimits_{\pm}
\REPEAT{
\left( \lm(1+h)\xi_1^{\pm}\cdot \lm(1+h^{-1})\xi_2^{\pm}\right)
}
% \tp\left(\lm(1+h)\eta_1^{\pm}\cdot \lm(1+h^{-1})\eta_2^{\pm}\right)
\;\frac{dh}{|h|}
\end{aligned}
\right.
\end{aligned}
\]
Comparing this with the combination of \eqref{eq:four pieces}, \eqref{eq:pm-pm} and \eqref{eq:pm-mp},
we have shown that 
$P_{\rm diag}(T_1)\boxtimes P_{\rm diag}(T_2) = P_{\rm diag}(T_1\boxtimes T_2)$ as required.
\end{proof}%end proof of prop.

\begin{dfn}[The diagonal subalgebra]
We define $\cA_{\rm diag}\defeq P_{\rm diag}(\cA)$. Note that by Proposition \ref{prop:Adiag-direct}, $\cA_{\rm diag}$ is a subalgebra.
\end{dfn}

We now examine the image of $\cA_{\rm diag}$ under the map~$\Psi:\cA\to \FA(\affR)$.
%%%%%%%%%%%%%%%%%%%%%%%%%%%%%%%%%%%%%

\begin{lem}\label{l:two halves}
\begin{romnum}
\item\label{li:TH1}
$\pi(0,-1)\xi(t)=\xi(-t)$. In particular, $\pi(0,-1)$ interchanges $\sH^+$ and $\sH^{-}$.
\item\label{li:TH2}
If $(b,a)\in \affRpos$ then $\pi(b,a)(\sH^\pm)\subseteq \sH^\pm$.
\item\label{li:TH3}
If $\xi,\eta\in \sH^{\pm}$ (i.e.\ both have the same ``parity'') then $\xi*_\pi \eta$ vanishes outside $\affRpos$.
\item\label{li:TH4}
If $\xi\in\sH^\pm$ and $\eta\in\sH^{\mp}$ (i.e.\ they have different ``parity'') then $\xi*_\pi\eta$ vanishes on $\affRpos$.
\end{romnum}
\end{lem}

The claims in the lemma follow easily from the definitions of $\pi$ and $\sH^{\pm}$, so we leave the details to the reader.

\begin{prop}[Intertwining projections]
\label{p:intertwine projections}
$\Psi P_{\rm diag} = P_e \Psi$ as maps $\sH\ptp\sHbar\to\FA(\affR)$.
\end{prop}

\begin{proof}
By linearity and continuity, it suffices to verify this identity on rank-one tensors in $\sH\tp\sHbar$. Let $\xi,\eta\in\sH$; then
\[ \begin{aligned}
\Psi(\xi\tp\eta)
& =\Psi\left(
  P^+\xi \tp P^+\eta
+ P^+\xi \tp P^-\eta
+ P^-\xi \tp P^+\eta
+ P^-\xi \tp P^-\eta
\right) \\
& =
  (P^+\xi)*_\pi(P^+\eta)
+ (P^+\xi)*_\pi(P^-\eta)
+ (P^-\xi)*_\pi(P^+\eta)
+ (P^-\xi)*_\pi(P^-\eta).
\end{aligned} \]
Hence by Lemma~\ref{l:two halves}\ref{li:TH4}, for every $(b,a)\in\affRpos$,
\[
P_e\Psi(\xi\tp\eta)(b,a)
= 
  (P^+\xi)*_\pi(P^+\eta)(b,a)
+ (P^-\xi)*_\pi(P^-\eta)(b,a)
= 
\Psi P_{\rm diag}(\xi\tp\eta)(b,a).
\]
Thus, $P_e\Psi(\xi\tp\eta)$ and $\Psi P_{\rm diag}(\xi\tp\eta)$ agree on $\affRpos$.

By definition, $P_e\Psi(\xi\tp\eta)$ vanishes outside $\affRpos$;
and by Lemma~\ref{l:two halves}\ref{li:TH4}, so does
$\Psi P_{\rm diag}(\xi\tp\eta)$.
We conclude that $P_e\Psi(\xi\tp\eta)= \Psi P_{\rm diag}(\xi\tp\eta)$ as required.
\end{proof}

Since $P_e:\FA(\affR)\to\FA(\affR)$ is a homomorphism with range $\FA_e(\affR)$
and $\Psi:\cA\to \FA(\affR)$ is an algebra isomorphism, this provides an alternative proof that $P_{\rm diag}$ is a homomorphism and $\cA_{\rm diag}$ is a subalgebra of~$\cA$.
Moreover, by the remarks at the start of this section, we may identify $\FA_e(\affR)$ with $\FA(\affRpos)$.
Thus $(\cA_{\rm diag}, \boxtimes)$ may be viewed as a realization of dual convolution for $\affRpos$.

\begin{rem}
Lemma~\ref{l:two halves}\ref{li:TH2} shows that the restriction of $\pi$ to $\affRpos$ splits as $\pi_{+}\oplus\pi_{-}$ where $\pi_{\pm}: \affRpos \to \cU(\sH^{\pm})$. Up to unitary equivalence, $\pi_+$ and $\pi_{-}$ are the only two infinite-dimensional unitary representations of $\affRpos$; they can also be constructed directly as induced representations.
Attempting to construct dual convolution for $\affRpos$ directly requires consideration of $\pi_+\tp \pi_+$, $\pi_+\tp\pi_-$, $\pi_-\tp\pi_+$ and $\pi_-\tp\pi_-$, and the fusion rules for the ``mixed parity'' cases are not so straightforward. Indeed, $\pi_+\tp\pi_-$ is not quasi-equivalent to an irreducible representation of $\affRpos$.
\end{rem}

Finally, we consider derivations on $\cA_{\rm diag}$ and hence on $\FA(\affRpos)$.
Let $D: \cA\to\cA^*$ be the derivation constructed in Section~\ref{s:derivation}.
Composing with the inclusion $\iota: \cA_{\rm diag}\to \cA$ and the restriction $\iota^*: \cA^*\to (\cA_{\rm diag})^*$, we obtain a derivation $D_1=\iota^*D\iota : \cA_{\rm diag} \to (\cA_{\rm diag})^*$.
Cyclicity and weak compactness of $D$ are inherited by $D_1$, just from the definition.
Now let $\top$ be the transpose operator from Definition~\ref{d:transpose}.
Since $\iota$ is a complete isometry and $\top_* \iota = \iota \top_*$,
complete boundedness of $\top D: \cA\to\cA^*$ implies complete boundedness of
$\top  D_1 =\iota^* \top  D \iota: \cA_{\rm diag} \to (\cA_{\rm diag})^*$.

It remains to check that $D_1$ is not identically zero.
Recall that by definition
\[
D(\xi_1\tp\eta_1)(\xi_0 \tp \eta_0)
=
\ip{S\xi_1}{ R\overline{\xi_0}}\ip{R\overline{\eta_0}}{\eta_1}
\]
where $(S\xi)(t)=\sgn(t)\xi(t)$ and $R\xi(t)=\xi(-t)$.
Fix some non-zero vector $\alpha\in \sH^+$ and put $\beta=\overline{R\alpha} \in\sH^{-}$, so that $\alpha\tp\alpha$ and $\beta\tp\beta$ belong to $\cA_{\rm diag}$. Since $S\alpha=\alpha$ and $R^2={\rm id}$ we have
\[
D(\alpha \tp \alpha)(\beta\tp\beta)
=
\ip{S\alpha}{R\overline{\beta}} \ip{R\overline{\beta}}{\alpha}
=
\ip{\alpha}{\alpha}\ip{\alpha}{\alpha}
\neq 0.
\]

Intertwining $D_1$ with $\Psi$ yields a non-zero derivation $\widetilde{D}: \FA_e(\affR)\to \FA_e(\affR)^*$.
This derivation turns out to coincide, up to a scaling factor, with the derivation constructed in \cite{Choi-Ghandehari14}; the proof requires the orthogonality relations for $\pi_{\pm}$ or the Plancherel theorem for $\affRpos$.
Since $\Psi$ is a completely isometric algebra isomorphism,
$\widetilde{D}$ inherits the properties of~$D_1$. In particular $\widetilde{D}$ is weakly compact and ``co-completely bounded'' (using the terminology of \cite{choi_Z2-pre}), properties which were less obvious from the original construction in \cite{Choi-Ghandehari14}.

\end{section}

%%%%%%%%%%%%%%%%%%%%%%%%%%%%%%%%%%%%%%%%%%%%%%%%%%

\begin{section}{A new Banach algebra structure on $L^p(\Rx)\ptp L^q(\Rx)$}
\label{s:Lp-version}

Throughout this section, we assume $1<p <\infty$ and denote by $q$ the conjugate index to~$p$.
We denote the usual pairing between $L^p(\Rx)$ and $L^q(\Rx)$ by $\pair{\quad}{\quad}_{p,q}$; note that $\pair{\xi}{\overline{\eta}}_{2,2} = \ip{\xi}{\eta}$.

For sake of precision, recall that there is an isometric, $\Cplx$-linear isomorphism of Banach spaces
\[
\tilde{\iota}: L^2(\Rx)\ptp \overline{L^2(\Rx)} \to L^2(\Rx)\ptp L^2(\Rx)
\qquad; \qquad  \xi\tp\eta \mapsto \xi\tp\overline{\eta}.
\]
By intertwining with $\tilde{\iota}$, we may transfer the Banach algebra structure defined on $L^2(\Rx)\ptp \overline{L^2(\Rx)}$ over to $L^2(\Rx)\ptp L^2(\Rx)$.
Moreover, one can use the natural analogue of the formula \eqref{eq:fusion of rank1} to equip $L^p(\Rx)\ptp L^q(\Rx)$ with a Banach algebra structure, in a way that extends the $p=q=2$ case.

That is: for $\xi,\xi'\in L^p(\Rx)$ and $\eta,\eta'\in L^q(\Rx)$, we claim that
\[
(\xi\tp\eta)\boxtimes (\xi'\tp\eta') \defeq
\int_{\Rx} (\lm(1+h)\xi \cdot \lm(1+h^{-1})\xi' )
      \tp (\lm(1+h)\eta\cdot \lm(1+h^{-1})\eta')
\,\frac{dh}{|h|}
\]
is a well-defined Bochner integral taking values in $L^p(\Rx)\ptp L^q(\Rx)$, and that $\boxtimes$ extends to a bounded bilinear map
\[
\left(
L^p(\Rx)\ptp L^q(\Rx)
\right)
\times
\left(
L^p(\Rx)\ptp L^q(\Rx)
\right)
\to
L^p(\Rx)\ptp L^q(\Rx)
\]
which is commutative and associative.
We denote the resulting commutative Banach algebra $(L^p(\Rx)\ptp  L^q(\Rx), \boxtimes)$ by $\cA_p$.

Most of the steps needed to justify this claim consist of routine modifications of the arguments in Section~\ref{ss:define dual conv}, so we shall not give full details here. We highlight some of the relevant technical points.
% In showing that $\boxtimes$ is well-defined, there are two key technical points to consider.
\begin{enumerate}
\renewcommand{\theenumi}{(S\arabic{enumi})}
\renewcommand{\labelenumi}{(S\arabic{enumi})}
\item\label{li:key Lp}
There is an $L^p$-analogue of Lemma~\ref{l:fund unitary}, with an isometry  $V_p:L^p(\Rx\times\Rx) \to L^p(\Rx\times\Rx)$ defined by the formula
\[
V_p(X)(h,t) = X\left( \frac{t}{1+h} , \frac{t}{1+h^{-1}}\right).
\]
As before, this shows that for $\xi,\xi'\in L^p(\Rx)$ the function
\[
F:h\mapsto \lm(1+h)\xi\cdot \lm(1+h^{-1})\xi'\]
is a.e.\ equal to a strongly measurable $L^p(\Rx)$-valued function, with
\[
\left(\int_{\Rx} \norm{F(h)}_p^p \,\frac{dh}{|h|} \right)^{1/p} = \norm{\xi}_p \norm{\xi'}_p \;.
\]
One then performs the same construction with $p$ replaced by $q$.
\item
However, one has to be careful taking pointwise products of two functions in $L^p(\Rx)$ or $L^q(\Rx)$. The expression defining $F(h)$ {\it a priori} only takes values in $L^{p/2}(\Rx)$, which for $1<p<2$ is not a Banach space (it is complete and quasi-normed, but not locally convex).
\item
Once one has shown that $\boxtimes$ is well-defined and contractive as a bilinear map $(L^p(\Rx)\ptp L^q(\Rx)) \times (L^p(\Rx)\ptp L^q(\Rx)) \times(L^p(\Rx)\ptp L^q(\Rx))$, one can prove it is commutative and associative by repeating the arguments of Section~\ref{ss:basic properties} almost verbatim; the key point is that $C_c(\Rx)$ is still norm-dense in both $L^p(\Rx)$ and $L^q(\Rx)$.
\item
The abstract description of $\boxtimes$ for $\sH\ptp\sHbar$, shown in Figure~\ref{fig:fancy POV}, has a natural and straightforward generalization to the $L^p$-setting, which is sketched in Figure~\ref{fig:Lp fancy POV}.
\end{enumerate}

\begin{figure}[htp]
\[
\begin{aligned}
& & (L^p(\Rx) \ptp L^q(\Rx) ) \ptp (L^p(\Rx)\ptp L^q(\Rx))
\\
&
\strut \xrightarrow{\rm shuffle} \strut &
(L^p(\Rx) \ptp L^p(\Rx) ) \ptp (L^q(\Rx)\ptp L^q(\Rx))
 \\
&
\strut \xrightarrow{\rm embed} \strut  &
L^p(\Rx\times \Rx)  \ptp L^q(\Rx\times \Rx)
 \\
&
\strut \xrightarrow{V_p \tp V_q} \strut  &
L^p(\Rx\times \Rx)  \ptp L^q(\Rx\times \Rx)
 \\
&
\strut \xrightarrow{\rm identify} \strut  &
L^p(\Rx; L^p(\Rx))  \ptp L^q(\Rx ; L^q(\Rx))
\\
&
\strut \xrightarrow{\rm diagonal} \strut  &
L^1(\Rx; L^p(\Rx)  \ptp L^q(\Rx))
\\
&
\strut \xrightarrow{\rm trace} \strut  &
L^p(\Rx)  \ptp L^q(\Rx)
\end{aligned}
\]
\caption{The $L^p$-analogue of Figure~\ref{fig:fancy POV}}
\label{fig:Lp fancy POV}
\end{figure}

For $p=2$,  $\tilde{\iota}:\cA\to \cA_2$ is an isometric isomorphism of Banach algebras.
Since $\Psi:\cA\to C_0(\affR)$ is an injective homomorphism, it follows that $\cA_2\cong\cA$ is semisimple, and that we can identify $\cA_2$ with a Banach function algebra on~$\affR$.
We now show that the same is true for~$\cA_p$.

The formula $\pi(b,a)\xi(t)\defeq e^{2\pi i bt} \xi(ta)$ still defines an isometric, SOT-continuous representation of $\affR$ on $L^p(\Rx)$.
Hence, for each $\xi\in L^p(\Rx)$ and $\eta\in L^q(\Rx)$ there is an associated coefficient function:
\begin{equation}
\Psi_p(\xi\otimes\eta)(b,a)\defeq \pair{\pi(b,a)\xi}{\eta}_{p,q}
=
\int_{\Rx} e^{2\pi ibt} \xi(ta)\eta(t)\,\frac{dt}{|t|}
\qquad((b,a)\in\affR).
\end{equation}
This formula defines a contractive linear map $\Psi_p: L^p(\Rx)\ptp L^q(\Rx) \to C_b(\affR)$.
(Note that $\Psi_2\circ\tilde{\iota}=\Psi$.)
We define $\FA^p_\pi$ to be the space $\Psi_p(\cA_p)\subset C_b(\affR)$ equipped with the quotient norm pushed forwards from $\cA_p/\ker\Psi_p$.
One can show that $\Psi_p:\cA_p\to C_b(\affR)$ is an algebra homomorphism, by a direct calculation using the $L^p$-analogue of Proposition~\ref{p:fusion of coeff}.
Hence $\FA^p_\pi$  is a Banach function algebra, which in the case $p=2$ is just $\FA_\pi(\affR)$.

\begin{prop}
$\FA^p_\pi \subseteq C_0(\affR)$.
\end{prop}

\begin{proof}
By linearity and continuity it suffices to prove that for each $\xi,\eta\in C_c(\Rx)$ the coefficient function $f\defeq \Psi_p(\xi\tp\eta)$ belongs to $C_0(\affR)$.
This now follows because $\xi,\eta\in L^2(\Rx)$ and $\FA^2_\pi = \FA_\pi(\affR)=\FA(\affR)\subset C_0(\affR)$.

An alternative argument, which does not rely on the equality $\FA_\pi(\affR)=\FA(\affR)$, goes as follows.
Since $\xi,\eta\in C_c(\Rx)$, there is a compact $K\subset\Rx$ such that $f$ is supported inside $\Real\times K$. Also, for each $a\in K$ we have $f(\underline{\quad},a)\in C_0(\Real)$, since $t\mapsto |t|^{-1} \xi(ta)\eta(t)$ is integrable (use the Riemann--Lebesgue lemma for the Fourier transform on $\Real$).
By a standard compactness argument, whose details we omit, we conclude that $f\in C_0(\Real\times K) \subset C_0(\affR)$.
\end{proof}

So far everything has been a straightforward translation of what was done for the $p=2$ case.
In contrast, the next result seems to require extra work.

\begin{thm}\label{t:Psi_p is injective}
$\Psi_p:L^p(\Rx)\ptp L^q(\Rx) \to C_0(\affR)$ is injective. Consequently, $\Psi_p:\cA_p\to\FA^p_\pi$ is an isometric isomorphism of Banach algebras.
\end{thm}

For $p=2$ this is a special case of general results already mentioned in Section~\ref{s:prelim}. For general $p$, we make use of results from \cite{eymard-terp} that are particular to $\pi$ and $\affR$.
Consider the following space:
\begin{equation}
\sV_0 \defeq \{ \xi \in \FA(\Real) \colon \operatorname{supp}\xi \text{ is compact and disjoint from $\{0\}$} \}.
\end{equation}
$\sV_0$ is a linear subspace of $\FA(\Real)$; standard properties of $\FA(\Real)$ imply that $\sV_0$ is norm-dense in $L^p(\Rx)$ for every $p\in (1,\infty)$.
The following lemma is a special case\footnotemark\ of a result from \cite{eymard-terp}, restated in a more direct form to avoid possible clashes of notational conventions.
\footnotetext{The lemma implies that for $\xi$ and $\eta$ in $\sV_0$, the corresponding rank-one operator belongs to $\pi(L^1(\affR))\subset \Bdd(L^2(\Rx))$, and this is the form in which Eymard and Terp state their result. In fact, they obtained a sharper result, which characterizes those $f\in L^1(\affR)$ such that $\pi(f)$ is a rank-one operator on~$L^2(\Rx)$.}

\begin{lem}[Eymard--Terp]\label{l:ET-rank-one}
Let $\xi,\eta\in \sV_0$. Then there exists $f\in L^1(\affR)$ such that, for every $\al,\beta \in C_c(\Rx)$,
\begin{equation}\label{eq:good vectors}
\int_{\affR} f(b,a) \left( \int_{\Rx} [\pi(b,a)\al ]\cdot \beta \right)  \, db\,\frac{da}{|a|^2}
=
\int_{\Rx} \al \eta \int_{\Rx} \xi \beta
\end{equation}
where the integrals over $\Rx$ are taken with respect to the Haar measure of this group.
\end{lem}

Since we only need a subset of Eymard and Terp's result, we include a proof of the lemma for the reader's convenience.

\begin{proof}[Proof (following {\cite[Prop.\ 1.13]{eymard-terp}})]
Let $\xi,\eta\in\sV_0$. The right-hand side of Equation~\eqref{eq:good vectors} is equal, after a change of variables $a\mapsto ta$, to
\begin{equation}\label{eq:internetdownagain}
\int_{\Rx}\int_{\Rx} \xi(t)\eta(ta) \ \alpha(ta) \beta(t)\, \frac{da}{|a|}\frac{dt}{|t|} \;.
\tag{$*$}
\end{equation}
If we can find $f\in L^1(\affR)$ such that $|a|^{-1}\int_{\Real} f(b,a)  e^{2\pi ibt}\,db= \xi(t)\eta(ta)$ for a.e.~$t,a\in\Rx$, then substituting this into \eqref{eq:internetdownagain} and using Fubini would give the left-hand side of Equation~\eqref{eq:good vectors}.

Let $g(t,a) = \xi(t)\eta(ta)$ viewed as a function $\Real\times\Rx\to\Cplx$.
The assumptions on $\xi$ and $\eta$ ensure that $a\mapsto g(\underline{\quad}, a)$ is a continuous function $\Rx\to\FA(\Real)$ which has compact support.
(C.f.~\cite[Exemple 1.17]{eymard-terp}).
Applying the inverse Fourier transform for $\Real$ to $g$ in the first variable, we obtain $f_1\in L^1(\Real\times \Rx, |a|^{-1}d(t,a))$ which satisfies
\[
\int_{\Real} f_1(b, a) e^{-2\pi ibt} \, db  = g(t,a)=\xi(t)\eta(ta).
\]
Thus the function $f(b,a) = |a| f_1(-b,a)$ has the required properties.
\end{proof}

\begin{proof}[Proof of Theorem~\ref{t:Psi_p is injective}]
It suffices to prove that $\Psi_p:L^p(\Rx)\ptp L^q(\Rx)\to C_0(\affR)$ is injective; the rest of the theorem follows from earlier observations.

Let $\xi,\eta\in \sV_0$ and let $f\in L^1(G)$ be as provided by Lemma~\ref{l:ET-rank-one}.
Let $j_p:L^p(\Rx)\ptp L^q(\Rx)\to \Bdd(L^p(\Rx))$ be the map which sends an elementary tensor $\al\tp\beta$ to the rank-one operator $\gamma \mapsto \pair{\gamma}{\beta}_{p,q}\alpha$. Then we may rewrite  Equation~\ref{eq:good vectors} as:
\[
\int_{\affR} f(b,a) \Psi_p(\al\tp\beta)(b,a)\,db\,\frac{da}{|a|^2}
=\pair{j_p(\al\tp\beta) \xi}{\eta}_{p,q}
\qquad\text{for all $\al,\beta\in C_c(\Rx)$.}
\]
Hence, by linearity and continuity of $j_p$ and $\Psi_p$,
\[
\int_{\affR} f(b,a) \Psi_p(w)(b,a)\,db\,\frac{da}{|a|^2}
=\pair{j_p(w) \xi}{\eta}_{p,q}
\qquad\text{for all $w\in L^p(\Rx)\ptp L^q(\Rx)$.}
\]

In particular, suppose $w\in \ker(\Psi_p)$. Then $\pair{j_p(w)\xi}{\eta}=0$. Since this holds for all $\xi,\eta\in \sV_0$, and since $\sV_0$ is norm-dense in $L^p(\Rx)$ and in $L^q(\Rx)$, it follows that $j_p(w)=0$. Since $L^p(\Rx)$ has the approximation property, $j_p$ is injective, and we conclude that $w=0$ as required.
\end{proof}

\begin{rem}\label{r:justify}
Define $\Phi_p: M(\affR) \to \Bdd(L^p(\Rx))$ defined as follows: for $\mu\in M(\affR)$, $\xi\in L^p(\Rx)$, and $\eta\in L^q(\Rx)$, let
\begin{equation}
\pair{\Phi_p(\mu) \xi}{\eta}_{p,q}
\defeq \int_{\affR} \pair{\pi(b,a) \xi}{\eta}_{p,q} \,d\mu(b,a).
\end{equation}
$\Phi_p$ is a contractive, weak$^*$-weak$^*$ continuous algebra homomorphism,
and it can be identified with the adjoint of~$\Psi_p$.
Hence injectivity of $\Psi_p$ is equivalent to weak$^*$-density of $\Phi_p(M(\affR))$ in $\Bdd(L^p(\Rx))$.
In effect, our proof of Theorem~\ref{t:Psi_p is injective} works by showing that $\Phi_p(L^1(\affR))$ contains a norm-dense subspace of ${\mathcal K}(L^p(\Rx))$ and hence is weak$^*$-dense in $\Bdd(L^p(\Rx))$. While this formulation of the proof is more intuitive, it does not seem to make the argument significantly simpler.
Note that for $p=2$ the weak$^*$-density result would follow from general facts about unitary representations of locally compact groups, but the proofs of those facts use ${\rm C}^*$-algebraic tools which are not available for general representations on $L^p$-spaces.
%%
%% The proof of this general result uses ${\rm C}^*$-algebraic tools such as Schur's lemma for unitary representations, the Kaplansky density theorem, and the von Neumann bicommutant theorem, which are not known for general representations on $L^p$-spaces.
\end{rem}

For any locally compact group $G$, the {\itshape Fig\`a-Talamanca--Herz algebra} $\FA_p(G)$ is defined to be the coefficient space of the left regular representation of $G$ on $L^p(G)$. Note that $\FA_2(G)=\FA(G)$.
%% $\FA_p(G)$ is usually regarded as the correct $L^p$-version of the Fourier algebra.
%%
We have seen above that $\FA^2_\pi=\FA_2(\affR)$; we now show that this fails for all other~$p$.

\begin{prop}\label{p:compare norms}
Let $p\in (1,2)\cup (2,\infty)$. There exists a sequence $(f_n)$ in $\FA^p_\pi \cap \FA_p(\affR)$ such that each $f_n$ has norm $1$ in $\FA^p_\pi$ but $f_n\to 0$ in $\FA_p(\affR)$.
\end{prop}

\begin{proof}
For this proof, we denote the norm in $\FA_p(\affR)$ by $\norm{\cdot}_{\FA_p}$.

Since $\affR$ is amenable, a result of Herz\footnotemark\ implies that $\FA_2(\affR)\subseteq \FA_p(\affR)$, and that the inclusion is norm-decreasing.
\footnotetext{For a guide to the relevant parts of Herz's papers, see the appendix of \cite{choi_dirfin15}. For a direct approach, see the proof of Theorem 8.3.9 in \cite{Der_book} and the historical notes which follow it.}%
Therefore, if we take $\xi \in (L^2\cap L^p)(\Rx)$ and $\eta\in (L^2\cap L^q)(\Rx)$ and set $f=\Psi_p(\xi\tp\eta)=\Psi_2(\xi\tp\eta)$,
we have
\[ \norm{f}_{\FA^p_\pi} = \norm{\xi\tp\eta}_{\cA_p} = \norm{\xi}_p \norm{\eta}_q
\;\text{and}\;
\norm{f}_{\FA_p} \leq \norm{f}_{\FA_2} = \norm{f}_{\FA^2_\pi} = \norm{\xi\tp\eta}_{\cA_2} = \norm{\xi}_2\norm{\eta}_2 \;.
\]

It therefore suffices to find functions $\xi_n$ and $\eta_n$ that lie in every $L^p(\Rx)$ and satisfy
\[
\frac{\norm{\xi_n}_2}{\norm{\xi_n}_p} \frac{\norm{\eta_n}_2}{\norm{\eta_n}_q}
 \to 0 \quad\text{as $n\to\infty$}\;,
\]
since we may then take $f_n = \norm{\xi_n}_p^{-1}\norm{\eta_n}_q^{-1}\Psi_p(\xi_n\tp\eta_n)$. 
Consider $\gamma_n:\Rx\to \{0,1\}$ defined by $\gamma_n = 1_{[e^{-n},e^n]}$.
Then $\gamma_n \in L^p(\Rx)$ for all $p\in (1,\infty)$ with $\norm{\gamma_n}_p^p = 2n$.
For $1<p<2$, taking $\xi_n=\gamma_n$ and $\eta_n = \gamma_1$ yields
\[
\frac{\norm{\xi_n}_2}{\norm{\xi_n}_p} \frac{\norm{\eta_n}_2}{\norm{\eta_n}_q}
= (2n)^{\frac{1}{2}-\frac{1}{p}} 2^{\frac{1}{2}-\frac{1}{q}} =n^{\frac{1}{2}-\frac{1}{p}}\to 0 \;;
\]
while for $2<p<\infty$, taking $\xi_n=\gamma_1$ and $\eta_n = \gamma_n$ yields
\[
\frac{\norm{\xi_n}_2}{\norm{\xi_n}_p} \frac{\norm{\eta_n}_2}{\norm{\eta_n}_q}
=  2^{\frac{1}{2}-\frac{1}{p}} (2n)^{\frac{1}{2}-\frac{1}{q}} = n^{\frac{1}{2}-\frac{1}{q}} \to 0 \;.
\]
So in both cases we have the desired sequences.
\end{proof}

\begin{thm}\label{t:not A_p}
If $p\in (1,2)\cup (2,\infty)$, then $\FA_p(\affR)\not\subseteq \FA^\pi_p$.
\end{thm}

\begin{proof}
Suppose that $\FA_p(\affR)\subseteq\FA^p_\pi$. Since both Banach spaces embed continuously in $C_0(\affR)$, the closed graph theorem would then imply that the inclusion map $\FA_p(\affR)\to \FA^p_\pi$ is continuous. But this contradicts Proposition~\ref{p:compare norms}.
\end{proof}

\begin{rem}
Since $\FA^p_\pi$ is the coefficient space of an isometric group representation on an $L^p$-space, it is contained in the multiplier algebra of $\FA_p(\affR)$. This follows from an $L^p$-version of Fell's absorption principle (valid for any locally compact group), which appears to be folklore and goes back to the 1960s/70s.
It would be interesting to study the relationship between $\FA^p_\pi$ and $\FA_p(\affR)$ in greater detail.
\end{rem}

\end{section}

\begin{section}{Concluding remarks}
\label{s:conclusion}
We finish by suggesting some avenues for further exploration.

\paragraph{Affine groups of other local fields.}
Much of \cite{eymard-terp} works in the general setting of a field ${\sf k}$ which is locally compact, second-countable and non-discrete, together with the corresponding affine group ${\sf k}\rtimes {\sf k}^\times$. All the calculations of Section~\ref{s:fusion-DC} and Section~\ref{s:Lp-version} should remain true for such a ${\sf k}$, provided that one replaces the exponential function in the definition of $\pi$ with a nontrivial character of $({\sf k},+)$.
However, Sections~\ref{s:derivation} and \ref{s:subalgebra} use certain special features of $\Rx$ that are not shared by ${\sf k}^\times$, and we do not expect them to generalize to ${\mathbb Q}_p$, for instance.

\paragraph{Constructing explicit derivations on Fourier algebras.}
The question of which groups $G$ allow non-zero derivations $\FA(G)\to\FA(G)^*$ has been intensively studied in recent years. The calculations in Section~\ref{s:derivation} may give new ideas or techniques for constructing derivations on Fourier algebras of other (Lie) groups.

\paragraph{A concrete model for LCQG questions.}
By enhancing the decomposition in Figure~\ref{fig:fancy POV} with operator-space structure, using row and column Hilbert spaces in the appropriate places, one can show that $\boxtimes$ extends to a completely contractive map $\cS_1(\sH) \ptp_{\rm op} \cS_1(\sH)\to\cS_1(\sH)$, where $\ptp_{\rm op}$ denotes the projective tensor product of operator spaces. The adjoint of this map is a $*$-homomorphism $\Delta_\boxtimes:\Bdd(\sH)\to \Bdd(\sH\tp^2\sH)$, which is coassociative since $\boxtimes$ is associative.
Moreover, the adjoint of $\Psi:\cS_1(\sH) \to \FA(\affR)$ coincides with the canonical $*$-homo\-morphism $\VN(\affR) \to \Bdd(\sH)$ obtained by sending $\lambda(b,a)$ to $\pi(b,a)$. Because $\Psi$ is a homomorphism, $\Psi^*$ intertwines $\Delta_\boxtimes$ with the canonical comultiplication $\Delta$ on $\VN(\affR)$.

It might be interesting to study various general constructions for Hopf von Neumann algebras,  using $(\Bdd(\sH),\Delta_{\boxtimes})$ as our concrete model of $(\VN(\affR),\Delta)$.
In particular, to our knowledge it remains an open question if the operator systems $\operatorname{WAP}(\widehat{G})$ and $\operatorname{LUC}(\widehat{G})$ are subalgebras of $\VN(G)$ for non-abelian $G$; our concrete model may provide a new angle of attack when $G=\affR$.

One note of warning: the transpose operator $\top: \Bdd(\sH)\to\Bdd(\sH)$ is \emph{not} intertwined with the canonical involution on $\VN(\affR)$, because $\top(\pi(b,a)) \neq \pi(b,a)^*$. If we wish to also introduce Kac algebra structure on $(\Bdd(\sH),\Delta_{\boxtimes})$, the antipode is given not by $\top$ but by a unitarily similar operator.

\paragraph{Questions regarding $\cA_p$ and $\FA^p_\pi$.}
Let $p\in (1,2)\cup (2,\infty)$. 
\begin{enumerate}
\renewcommand{\labelenumi}{Q\arabic{enumi}.}
\item
Does $\cA_p$ have a bounded approximate identity?
\item
Is $\cA_p$ weakly amenable? In Section~\ref{s:derivation} we wrote down an explicit $\Phi\in(\cA\ptp\cA)^*$ which defines a non-zero derivation $\cA\to\cA^*$. However, $\Phi$ does not extend to a bounded bilinear functional on $\cA_p\ptp\cA_p$.
\item
Is $\FA^p_\pi$ \emph{natural} as a Banach function algebra on $\affR$? Equivalently: are all characters on $\cA_p$ of the form $w \mapsto \Psi_p(w)(b,a)$ for some $(b,a)\in\affR$?
\item
Assuming a positive answer to the previous question: which other function-algebra properties of $\FA(\affR)$ are shared by $\FA^p_\pi$? For example: is this algebra regular? Tauberian? 
It is not clear to the authors if $\FA^p_\pi$ contains any non-zero compactly supported functions.
\end{enumerate}

\subsection*{Acknowledgments}
This work was initiated during a visit of Choi to the University of Delaware in September 2017, supported by LMS grant 41669 from the London Mathematical Society (Scheme 4), and continued during a second visit in September 2018, supported by the Visitor Fund of the Department of Mathematics and Statistics at Lancaster University.
Ghandehari was supported by  NSF grant DMS-1902301 while this work was being completed.

\end{section}

\appendix

\begin{section}{Tensor products of induced representations}\label{app:induced}
Consider a semidirect product $G=N \rtimes H$. The left action of $H$ on $N$ is denoted by $h\cdot n$; if $\sigma \in \widehat{N}$ then the corresponding \emph{left} action of $H$ on $\widehat{N}$ is defined by $h\cdot \sigma :n \mapsto \sigma(h^{-1}\cdot n)$.

Given a (continuous, unitary) representation $\sigma: N \to \cU(\cH_\sigma)$ we \emph{define} the induced representation $\Ind_N^G\sigma : G \to \cU(L^2(H,\cH_\sigma))$ by the formula
\begin{equation}\label{def:induced}
\Ind_N^G\sigma(n,h)\xi(k) = (k\cdot\sigma)(n) [\xi(h^{-1}k)] = \sigma(k^{-1}\cdot n) [ \xi(h^{-1}k)],
\end{equation}
for $\xi \in L^2(H,\cH_\sigma)$,  $n\in N$,  $h,k \in H$.
Combining Proposition 2.41 and Theorem 2.58 of \cite{KTbook}, we get the following theorem:
\begin{thm}\label{thm:tensor-formula}
Let $G=N\rtimes H$ and let $\pi_1$ and $\pi_2$ be representations of $N$. For $i=1,2$ let $\Pi_i=\Ind_{\pi_i}$ be the induced representation of $G$ on $L^2(H,\cH_{\pi_i})$. Then
\[
\Pi_1\otimes \Pi_2\simeq\int_H^\oplus\Ind_N^G(h\cdot\pi_1\otimes \pi_2)\,dh
\]
via the unitary map
$W:L^2(H,\cH_{\pi_1})\otimes^2L^2(H,\cH_{\pi_2})\rightarrow \int_H^\oplus L^2(H,\cH_{\pi_1}\otimes^2\cH_{\pi_2})$
defined by
\[
\ W(f\otimes g)=\int_H^\oplus\phi_{\rho(h)(f)\otimes g}\, dh,
\]%
where $\rho$ is the right regular representation ($\rho(h)f(k)=f(kh)$), and $\phi_{f\otimes g}\in L^2(H,\cH_{\pi_1}\otimes^2\cH_{\pi_2})$ is defined by 
$\phi_{f\otimes g}(h)=f(h)\otimes g(h)$.
\end{thm}
\begin{proof}
A direct calculation shows that $W$ preserves the inner product:
\[ \begin{aligned}
\langle W(f\otimes g), W(f'\otimes g')\rangle
&= \int_H\langle \phi_{\rho(h)f\otimes g}, \phi_{\rho(h)f'\otimes g'}\rangle_{L^2(H,\cH_{\pi_1}\otimes^2\cH_{\pi_2})} \,dh\\
&= \int_H \int_H\langle \phi_{\rho(h)f\otimes g}(k), \phi_{\rho(h)f'\otimes g'}(k)\rangle_{\cH_{\pi_1}\otimes^2\cH_{\pi_2}} \,dk\, dh\\
&= \int_H \int_H \langle\rho(h)f(k)\otimes g(k), \rho(h)f'(k)\otimes g'(k)\rangle_{\cH_{\pi_1}\otimes^2\cH_{\pi_2}} \,dk \, dh\\
&= \int_H \int_H \langle\rho(h)f(k), \rho(h)f'(k)\rangle\langle g(k), g'(k)\rangle \,dk \, dh\\
&= \int_H \int_H \langle f(kh), f'(kh)\rangle dh\langle g(k), g'(k)\rangle \,dk\\
&=  \langle f, f'\rangle \langle g,g'\rangle.
\end{aligned} \]
Using \eqref{def:induced}, it is easy to verify that  $\rho(h)\Ind_N^G\pi_1(x) f=\Ind_N^G(h\cdot \pi_1)(x)(\rho(h)f)$. Thus,
\[ \begin{aligned}
W(\Pi_1\otimes\Pi_2)(x)(f\otimes g)&= W(\Pi_1(x)f\otimes \Pi_2(x)g)\\
&= \int_H^\oplus \phi_{\rho(h)(\Pi_1(x)f)\otimes \Pi_2(x)g} \, dh\\
&= \int_H^\oplus \phi_{\Ind_N^G(h\cdot \pi_1)(x)(\rho(h)f)\otimes \Pi_2(x)g} \,dh\\
&= \int_H^\oplus \Ind_N^G(h\cdot \pi_1\otimes \pi_2)(x)\phi_{\rho(h)f\otimes g} \,dh\\
&= \left(\int_H^\oplus \Ind_N^G(h\cdot \pi_1\otimes \pi_2)\, dh\right)(x)(\int_H^\oplus \phi_{\rho(h)f\otimes g} \,dh)\\
&= \left(\int_H^\oplus \Ind_N^G(h\cdot \pi_1\otimes \pi_2)\, dh\right)(x)(W(f\otimes g)).
\end{aligned} \]
\end{proof}

As an application, we now derive an alternative proof of Proposition~\ref{p:fusion of coeff}. We use the same notation as defined in Section~\ref{s:prelim}.
% We denote the unitary dual of $\Real$ by $\widehat{\Real}$.
%
For $r\in \Real$, define $\chi_r : \Real\to\Cplx$ by $\chi_r(t) = \exp(2\pi i rt)$, so that $r\mapsto \chi_r$ is a group isomorphism $\Real\to\widehat{\Real}$.
%
% This gives an explicit isomorphism $\Real\to \widehat{\Real}$, $r\mapsto \chi_r.$

\begin{cor}[Explicit fusion relation for $\pi$]
Let $\xi$, $\xi'$, $\eta$, $\eta'\in \sH$. Then
\begin{equation}
\begin{gathered}
\ip{\pi(b,a)\xi}{\eta}  \ \ip{\pi(b,a)\xi'}{\eta'} \\
= \int_{\Real^\times}   \ip{ \pi(b,a) ( \lambda(1+r)\xi\cdot \lambda(1+r^{-1})\xi')\ }{ \  \lambda(1+r)\eta\cdot \lambda(1+r^{-1})\eta' }  \;\frac{dr}{|r|},
\end{gathered}
\end{equation}
where $\lambda$ is the left regular representation of $\Rx$ on~$\sH$.
\end{cor}

\begin{proof}
As in Equation \eqref{eq:induced}, we denote $\Pi=\Ind_\Real^\affR\chi_1$. Recall that $\pi= J \Pi(\cdot) J$, where $Jf=\check{f}$. 
By Theorem \ref{thm:tensor-formula},
\[ W(\Pi\otimes \Pi)(b,a)W^{-1}=\int_{\Real^\times}^\oplus\Ind_\Real^\affR \chi_{\frac{1}{r}+1}(b,a)\, \frac{dr}{|r|}
\]
via the unitary map
\[
W:L^2(\Real^\times)\otimes^2L^2(\Real^\times)\rightarrow \int_{L^2(\Real^\times)}^\oplus L^2(\Real^\times)\, \frac{dr}{|r|},
\quad
W(f\otimes g)=\int_{\Real^\times}^\oplus (\rho(r)f) g\, \frac{dr}{|r|},
\]
where $\rho$ is the right regular representation of~$\Rx$, i.e. $\rho(r)f(s)=f(sr)$. Here, we have used the fact that $r\cdot\chi_1\otimes\chi_1\simeq\chi_{\frac{1}{r}+1}$. Hence
\[
\begin{aligned}
\ip{\pi(b,a)\xi}{\eta}  \ \ip{\pi(b,a)\xi'}{\eta'}
& = \ip{\Pi(b,a)J\xi}{J\eta}  \ \ip{\Pi(b,a)J\xi'}{J\eta'} \\
& = \ip{ (\Pi\otimes\Pi) (b,a)(\check\xi\otimes\check{\xi'})}{ \check\eta\otimes\check{\eta'} } \\
\end{aligned}
\]
which expands out to
\begin{eqnarray*}
& \phantom{=} & \left\langle\left(\int_{\Real^\times}^\oplus\Ind_\Real^\affR \chi_{\frac{1}{r}+1}(b,a)\, \frac{dr}{|r|}\right)W(\check\xi\otimes\check{\xi'}),W(\check\eta\otimes\check{\eta'})\right\rangle \\
&=&\int_{\Real^\times}   \left\langle\Ind_\Real^\affR \chi_{\frac{1}{r}+1}(b,a)\left((\rho(r)\check\xi)\check{\xi'}\right), (\rho(r)\check\eta)\check{\eta'}\right\rangle\, \frac{dr}{|r|}\ \\
&=&\int_{\Real^\times} \int_{\Real^\times}  \exp(\frac{2\pi i b(1+r)}{sr}) \check\xi(\frac{sr}{a})\check{\xi'}(\frac{s}{a}) \overline{\check\eta(sr)}\overline{\check{\eta'}(s)}\, \frac{ds}{|s|} \frac{dr}{|r|}\ \\
&=&\int_{\Real^\times}  \int_{\Real^\times}   \exp(\frac{2\pi i b(1+r)}{sr}) \xi(\frac{a}{sr}){\xi'}(\frac{a}{s}) \overline{\eta(\frac{1}{sr})}\overline{{\eta'}(\frac{1}{s})} \, \frac{ds}{|s|}\, \frac{dr}{|r|} \\
&=&\int_{\Real^\times}  \int_{\Real^\times}   \exp(2\pi i bs) \xi({\frac{sa}{1+r}}){\xi'}(\frac{sa}{1+r^{-1}}) \overline{\eta(\frac{s}{1+r})}\overline{{\eta'}(\frac{s}{1+r^{-1}})}\, \frac{ds}{|s|}\, \frac{dr}{|r|}\ \\
&=& \int_{\Real^\times}  \ip{ \pi(b,a) ( \lambda(1+r)\xi\cdot \lambda(1+r^{-1})\xi')\ }{ \  \lambda(1+r)\eta\cdot \lambda(1+r^{-1})\eta' }  \;\frac{dr}{|r|}
\end{eqnarray*}
as required.
\end{proof}

\end{section}

%\bibliography{istanbul-ref}
%\bibliography{ref}
%\bibliographystyle{alpha}

\vfill

\newcommand{\address}[1]{{\small\sc#1.}}
\newcommand{\email}[1]{\texttt{#1}}

\noindent
\address{Yemon Choi,
Department of Mathematics and Statistics,
Lancaster University,
Lancaster LA1 4YF, United Kingdom} 

\email{y.choi1@lancaster.ac.uk}

\noindent
\address{Mahya Ghandehari,
Department of Mathematical Sciences,
University of Delaware,
Newark, Delaware 19716, United States of America}

\email{mahya@udel.edu}

\end{document}